\documentclass[11pt,a4paper]{article}
\usepackage[latin1]{inputenc}
\usepackage{amsmath}
\usepackage{graphicx}
\usepackage{multirow}
\usepackage{mathrsfs}
\usepackage{amsfonts}
\usepackage{variations}
\usepackage{dsfont}
\usepackage{amssymb}
\author{BERNARD Damien}
\title{\textbf{Modular case of Levinson's theorem}}
\date{}
\def\C{\mathbb{C}}
\def\R{\mathbb{R}}

\def\a{\mathfrak{a}}
\def\b{\mathfrak{b}}

\def\l{\lambda_f}

\newtheorem{proposition}{Proposition}

\newtheorem{corollary}{Corollary}
\newtheorem{theorem}{Theorem}
\newcounter{remark}
\newenvironment{remark}{\addtocounter{remark}{1} \textbf{Remark \theremark}}{}
\newtheorem{lemma}{Lemma}

\newtheorem{conjecture}{Conjecture}
\begin{document}
\maketitle

\begin{abstract} We evaluate the integral mollified second moment of $L$-functions of primitive cusp forms and we obtain, for such $L$-function, an explicit positive proportion of zeros which lie on the critical line.
\end{abstract}

\tableofcontents

\section{Introduction and overview of the results}
 
Nowadays, we know that more than $41 \%$ of the non-trivial zeros of the Riemann zeta-function lie on the critical line (\cite{BCY} and \cite{Fe}). This theorem is the best one of a sequence of results about the percentage of zeros $\rho$ satisfying $\Re \rho = 1/2$. 

Historically, Selberg (\cite{S}) was the first one to show that this proportion is not zero without quantifying it. According to Titchmarsh (\cite{Ti}, part 10.9), it was calculated later on in Min's dissertation, that the proportion obtained by Selberg's method is very small. One may refer to the introduction of \cite{Ste} (page 8) for numerical values. In 1974, Levinson (\cite{L}) succeeded in proving that at least one third of the non-trivial zeros lie on the critical line by perturbing the Riemann zeta function by a linear combination of its derivatives.  A significant improvement, due to Conrey (\cite{C}), increased this proportion up to more than two fifths. In order to do this, he improves the general result of \cite{BCHB} on the asymptotic behaviour of the mollified second moment of the Riemann zeta function, when the coefficients of the mollifier are essentially given by the Moebius function, which allows him to work with a longer mollifier than Levinson's one. From this last result of Conrey and using a two parts mollifier, Bui, Conrey and Young (\cite{BCY}) proved that $41 \%$ of non-trivial zeros $\rho$ satisfies  $\Re \rho = 1/2$. 

Since the Riemann zeta-function is a $L$-function of degree one, it is rather natural to generalise these results to $L$-functions of higher degrees. For instance, Hafner has extended Selberg's result to $L$-functions of degree two. Precisely, if $f$ is a holomorphic cusp form of even weight and full level or an even Maa$\mbox{\ss}$  form of full level, let   $N_f(T)$ (resp. $N_{f,0} (T)$) be the number of the non-trivial zeros $\rho$ (resp. on the critical line) of $L(f,s)$ with $0<\Im (\rho)\leq T$.  Then, Hafner proved in \cite{H} and \cite{H2} that there exists a positive number $A$ such that  $N_{f,0} (T) > A N_f (T)$  for large $T$. Rezvyakova (\cite{Rez}) adapted (\cite{H}) to $L$-functions attached to automorphic cusp forms for congruence subgroups. Nevertheless, they do not give any explicit value for $A$ but by analogy with the Riemann zeta-function case, this constant should probably be close to zero. 

In \cite{Fa1}, Farmer applied Levinson's method to $L$-functions of an holomorphic cusp form $f$ of even weight and full level and succeeded in 
determining the asymptotic behaviour of the mollified integral second moment of $L(f,s)$ when the mollifier is a Dirichlet polynomial of length less than $T^{1/6-\varepsilon}$. From this result, Farmer obtained explicit lower bounds for the proportion of simple zeros of the $j^{th}$ derivative ($j\geq 1$) of the completed $L$-function of $L(f,s)$, which are on the critical line. Unfortunately, the length of the mollifier is too small to exhibit an explicit positive proportion of simple zeros on the critical line for $L(f,s)$ itself. Nevertheless, even thought Farmer did not remark it, his result proves that at least $1,65 \%$ of zeros of $L(f,s)$ satisfy $\Re s=1/2$ (confer part \ref{sec:EFPZCL}).

In this paper, we exhibit a positive proportion of zeros, which lie on the critical line for $L$-functions of  holomorphic primitive cusp forms. To get this result, we study the asymptotic behaviour of the smooth mollified second moment of $L(f,s)$ following  the method developed in \cite{Yo}. We choose a mollifier $\psi$, which is a Dirichlet polynomial of length $T^\nu$ defined in relation (\ref{eqn:defmol}) page \pageref{eqn:defmol}. Moreover, we introduce a smooth function $w$ compactly supported in $[T/4,2T]$ with some conditions on its derivatives (confer relations (\ref{eqn:hypwa}), (\ref{eqn:hypwb}) and (\ref{eqn:hypwc})). We prove the following theorem.
\begin{theorem} \label{pro:momentdordre2ramollilisse} Let $f$ be a holomorphic primitive cusp form of even weight, square-free level and trivial character. If $0<\nu< \frac{1-2\theta}{4+2\theta}$ and if  $\alpha, \beta$ are complex numbers satisfying  $\alpha, \beta \ll L^{-1}$ with $|\alpha+\beta|\gg L^{-1}$, then
\begin{small}
$$\int_{-\infty}^{+\infty}w(t)L\left(f,\frac{1}{2}+\alpha+it\right)L\left(f,\frac{1}{2}+\beta-it\right)\left| \psi(\sigma_0+it) \right|^2dt = \widehat{w}(0)c(\alpha, \beta) +O(T(\ln L)^4/L)$$
\end{small}where
\begin{small}
\begin{eqnarray}
 c(\alpha,\beta)=1+ \frac{1}{\nu}\frac{1-T^{-2(\alpha+\beta)}}{(\alpha+\beta)\ln T} \left.\frac{d^2}{dxdy}\left[M^{-\beta x-\alpha y}\int_{0}^{1}P(x+u)P(y+u)du  \right]\right|_{x=y=0} \label{eqn:limiteM2}
 \end{eqnarray}
\end{small}and where $\theta=7/64$ refers to the exponent in the approximation towards Ramanujuan-Petersson-Selberg conjecture (confer relation \emph{(\ref{eqn:defRP})}).
\end{theorem}

\begin{corollary} \label{cor:proportion} Let $f$ be a holomorphic primitive cusp form of even weight, square-free level and trivial character. At least $ 2,97 \%$ of non-trivial zeros of $L(f,s)$ lie on the critical line $\Re (s)=1/2$. Assuming the Selberg conjecture, we improve  this percentage to $6,93\%$. In other words,
\begin{equation*}
\liminf_{T \rightarrow +\infty}\frac{N_{f,0}(T)}{N_f(T)} \geq  \begin{cases} 0,0297 & \text{unconditionally,} \\  0,0693 & \text{under Selberg conjecture. } \end{cases}
\end{equation*}

\end{corollary} 
\begin{remark}  As Heath-Brown (\cite{HB}) and Selberg pointed out, the study of the second mollified moment allows to obtain a lower bound for the proportion of the simple non-trivial zeros lying on the critical line. Unfortunately, in our case, the length of the mollifier is too small to get a positive proportion of simple zeros satisfying the Riemann hypothesis for now. We plan to get back to this issue in the close future.
\end{remark} 

\paragraph*{} The method we used can also be applied to determine the asymptotic behaviour of the smooth second moment of $L(f,s)$ close to the critical line. 
\begin{theorem} \label{th:moment2nonramolli}Let $f$ be a holomorphic primitive cusp form of even weight, square-free level $N$ and trivial character. If $\alpha, \beta$ are complex numbers satisfying  $\alpha, \beta \ll L^{-1}$, we have 
\begin{small}\begin{align*}\int_{-\infty}^{+\infty}w(t)L\left(f,\frac{1}{2}+\alpha+it\right)L\left(f,\frac{1}{2}+\beta-it\right) dt =  \a_f \int_{-\infty}^{+\infty} w(t)\ln t dt + \left[\b_f+\a_f \ln\left( \frac{\sqrt{N}}{2\pi}\right)\right] \widehat{w}(0) \\ +O\left(|\alpha+\beta|T(\ln T)^2+T^{\frac12+\theta+\varepsilon} \right). \end{align*}\end{small}with
\begin{small}
\begin{align*}
\a_f =\frac{12 N}{\pi^2\nu(N)}L\left(Sym^2f,1\right) \mbox{ \normalsize{and} }
\b_f = \frac{12 N}{\pi^2\nu(N)}L\left(Sym^2f,1\right) \left( \frac{L'\left(Sym^2f,1\right)}{L\left(Sym^2f,1\right)}+\gamma+\sum_{p|N}\frac{\ln p}{p+1}-\frac{2\zeta'(2)}{\zeta(2)}\right)
\end{align*} \end{small}and where $\nu(N) =  N\prod_{p|N}\left( 1+\frac{1}{p}\right)$. 
\end{theorem}
\begin{remark} Our result is non-trivial only in the case $|\alpha+\beta|=o(1/ \ln T)$, and furthermore, we need $|\alpha+\beta|=o(1/ \ln^2 T)$ to ensure that the term of order $T$ is significant. 
\end{remark}
\begin{corollary}\label{cor:moment2nonramollidc} Let $f$ be a holomorphic primitive cusp form of even weight, square-free level $N$ and trivial character. Then
\begin{small}\begin{align*}\int_{-\infty}^{+\infty}w(t)\left|L\left(f,\frac{1}{2}+it\right)\right|^2dt = \int_{-\infty}^{+\infty} w(t)\left[\a_f \ln\left( \frac{t\sqrt{N}}{2\pi}\right) +\b_f \right]dt +O\left(T^{\frac12+\theta+\varepsilon} \right). \end{align*}\end{small}
\end{corollary}
\begin{remark}
This corollary is in agreement with the conjecture stated in \cite{CFKRS} about integral moments of $L$-functions (see section \ref{ss:conjCFKRS}).
\end{remark}
\paragraph*{}In \cite{Z}, Zhang succeeded in determining the main term of this integral second moment of  $L(f,s)$ (without smooth function $w$) on the critical line.  Thanks to corollary \ref{cor:moment2nonramollidc}, we improve his result with the following more precise asymptotic expansion. When $f$ is a holomorphic cusp form of even weight for the full modular group, we can also refer to \cite{Go} where a similar asymptotic expansion is given.
\begin{corollary} \label{cor:corZhang}Let $f$ be a holomorphic primitive cusp form of even weight, square-free level $N$ and trivial character. We have
\begin{small}$$\int_{0}^{T}\left|L\left(f,\frac{1}{2}+it\right)\right|^2dt =  \a_f T\ln T + \left[\b_f+\a_f \ln\left( \frac{\sqrt{N}}{2\pi e}\right)\right]T +O\left(T/ \ln T \right). $$\end{small}
\end{corollary}

\paragraph*{Notations}
\begin{flushleft}
- If $f$ and $g$ are some functions of the real variable,  $f(x) \ll_A g(x) $ or $f(x)=\noindent O_A(g(x))$ mean that $|f (x)|$ is
smaller than a constant, which only depends on $A$, times $|g(x)|$ for large $x$.\\
- Similarly, the notation $f \asymp g$ means $f (x)\ll g(x)$ and $g(x) \ll f(x)$.\\
- The function $\ln$ refers to the natural logarithm function.\\
- From the Riemann zeta function $\zeta$, we define for any positive square-free integer $N$
\end{flushleft}
 $$\zeta^{(N)}(s)=\prod_{p|N}\left(1-\frac{1}{p^s}\right)\zeta(s).$$

\paragraph*{Acknowledgements} I would like to thank Guillaume Ricotta and Emmanuel Royer for all their comments and their remarks during this work. I want to express my gratitude to them for their kindness and for all their encouragements. I wish to thank the organisers of the Number Theory Seminar of Institut de Mathématiques de Bordeaux for their warm welcome in may 2012.

\section{Review on $L$-functions of primitive cusp forms}

For this section, we may refer to \cite{IK} (chapter $14$). Throughout this paper, $f$ denotes a holomorphic primitive cusp form of even weight $k$ and square-free level $N$. The Fourier expansion at the cusp
 $\infty$ of $f$ is given by
  $$f(z)=\sum_{n\geq 1} \lambda_f(n)n^{\frac{k-1}{2}}e^{2i\pi n z}$$
 for every complex number $z$ in the upper-half plane with the arithmetic normalisation $\lambda_f(1)=1$.
The Fourier coefficients $\lambda_f(n)$ satisfy the multiplicative relations 
\begin{eqnarray} 
 \lambda_f(n)\lambda_f(m)  & = &  \sum_{   \begin{array}{c} \scriptstyle{ d|(m,n)} 
 \\  \scriptstyle{(d,N)=1} \end{array} } \lambda_f \left( \frac{mn}{d^2} \right)  \label{eqn:multiplicativite0} \\ 
   \lambda_f(mn)  & = &  \sum_{   \begin{array}{c} \scriptstyle{d|(m,n) }
 \\  \scriptstyle{(d,N)=1} \end{array}} \mu(d)\lambda_f \left(\frac{m}{d}\right)\lambda_f \left(\frac{n}{d}\right) \label{eqn:multiplicativite} 
  \end{eqnarray}
 for all positive integers $m$ and $n$. Since $\lambda_f(1)\neq 0$, we may define the convolution inverse $(\mu_f(n))$ of the sequence $(\lambda_f(n))$. This is an arithmetic multiplicative function, which satisfies for a prime number $p$
\begin{eqnarray} \label{eqn:proprietemu1}
& & \mu_f(1)=1, \hspace{4mm}\mu_{f}(p)=-\lambda_f(p), \hspace{4mm}\mu_{f}(p^2)=\lambda_f(p)^2-\lambda_f(p^2)= \left\{ \begin{array}{ll}1 & \mbox{ if } p\nmid N \\ 0 & \mbox{ otherwise} \end{array} \right. \\
& & \mbox{and if  }j\geq3 \hspace{3mm}\mu_f(p^j)=0 . \label{eqn:proprietemu2}
\end{eqnarray}
We consider
$$L(f,s)=\sum_{n\geq 1}\frac{\lambda_f(n)}{n^s}= \prod_{p\in \mathcal{P}}\left(1-\frac{\lambda_f(p)}{p^{s}}+\chi_0 (p)\frac{1}{p^{2s}}  \right)^{-1}= \prod_{p\in \mathcal{P}}\left(1-\frac{\alpha_f(p)}{p^s} \right)^{-1} \left(1-\frac{\beta_f(p)}{p^s} \right)^{-1}  , $$
which is an absolutely convergent and non-vanishing Dirichlet series, an Euler product on $\Re (s)>1$, where $\chi_0$ denotes the trivial character of modulus $N$ and $\alpha_f(p)$, $\beta_f(p)$ are the complex roots of the quadratic equation $X^2-\lambda_f(p)X+\chi_0(p)=0$.  Moreover, the function 
 $$\Lambda(f,s)= \left( \frac{\sqrt{N}}{2\pi} \right)^s \Gamma\left( s+\frac{k-1}{2} \right)L(f,s)=L_{\infty}(f,s)L(f,s)$$  is the completed $L$-function of $L(f,s)$. It can be extended to an holomorphic function on $\C$ and satisfies the functional equation 
$$\Lambda(f,s)=\varepsilon(f)\Lambda(f,1-s)$$
where $\varepsilon(f)=\pm 1$. Remark that, by the duplication formula for the gamma function, the local factor at infinity can be written
\begin{eqnarray} \label{eqn:facteurlocalinfiniFM}
L_{\infty}(f,s) =\left(\frac{2^k}{8\pi}\right)^{1/2}\left(\frac{\sqrt{N}}{\pi}\right)^s \Gamma\left(\frac{s}{2}+\frac{k-1}{4} \right)\Gamma\left(\frac{s}{2}+\frac{k+1}{4} \right).
\end{eqnarray}

\section{Mollified second moment of $L$-functions of modular forms}

This section contains the proof of theorem \ref{pro:momentdordre2ramollilisse}. We define the sequence $(\mu_f(n))_{n\geq 1}$ as the convolution inverse of the sequence $(\lambda_f(n))_{n\geq 1}$ and we define a mollifier $\psi$ of the shape
\begin{eqnarray} \label{eqn:defmol} \psi(s)=\sum_{n\leq M}\frac{\mu_f(n)}{n^{s+\frac{1}{2}-\sigma_0}}P\left(\frac{\ln M/n}{\ln M} \right)
\end{eqnarray}
with $M=T^\nu$, $\sigma_0=\frac{1}{2}-\frac{R}{\ln T}$ where $R$ is a positive real number and $P$ is a real polynomial satisfying $P(0)=0$, $P(1)=1$. In addition, we choose a function $w: \R \mapsto \R$, which satisfies
 \begin{subequations} \label{eqn:hypw}
 \begin{align}
 & w \mbox{ is smooth,} \label{eqn:hypwa}\\
 & w \mbox{ is compactly supported with supp }w \subset [T/4;2T], \label{eqn:hypwb} \\
 & \mbox{for each }j\geq 0, \mbox{ we have }w^{(j)}(t) \ll_j \Delta^{-j}\mbox{ where }\Delta=T/L \mbox{ and }L=\ln T. \label{eqn:hypwc} 
 \end{align}
\end{subequations}
For more convenience, we set  
$$ I_f(\alpha,\beta)= \int_{-\infty}^{+\infty}w(t)L\left(f,\frac{1}{2}+\alpha+it\right)L\left(f,\frac{1}{2}+\beta-it\right)\left| \psi(\sigma_0+it) \right|^2dt . $$
To study the asymptotic behaviour of $I_f(\alpha,\beta)$, we need an explicit expression of $L(f,s)$ with $0\leq \Re(s)\leq 1$. In lemma \ref{lem:EFA}, we get an exact formula, also called ``approximate functional equation'', which gives an expression for $L(f,s+it)L(f,s-it)$ with $s$ in the critical strip where we can not use the Dirichlet series. Thanks to this new relation, we may split  $I_f(\alpha,\beta)$ in a diagonal term (without oscillation) and an off-diagonal term (with oscillation). The off-diagonal contribution is bounded in part \ref{subsec:ODT} whereas the  diagonal term is estimated in part \ref{subsec:EDT}.

\begin{lemma} \label{lem:EFA}
Let $G$ be any entire function, which decays rapidly in vertical strips, even and normalised by $G(0)=1$.  Then for each complex numbers $\alpha$, $\beta$ such that $0 \leq |\Re (\alpha) |, |\Re(\beta)| \leq 1/2$, we have
\begin{small}
\begin{eqnarray*}
& & L\left(f,\frac{1}{2}+\alpha+it\right)L\left(f,\frac{1}{2}+\beta-it\right)  =  \underset{m,n\geq 1}{\sum\sum} \frac{\lambda_{f}(m)\lambda_f(n)}{m^{\frac12+\alpha}n^{\frac12+\beta}} \left( \frac{m}{n}\right)^{-it}V_{\alpha,\beta}(mn,t)\\
    & &\hspace{68mm} +X_{\alpha,\beta,t}\underset{m,n\geq 1}{\sum\sum} \frac{\lambda_{f}(m)\lambda_f(n)}{m^{\frac12-\beta}n^{\frac12-\alpha}} \left( \frac{m}{n}\right)^{-it}V_{-\beta,-\alpha}\left(mn,t\right)
\end{eqnarray*}\end{small}where  
\begin{small}
\begin{eqnarray*} \label{eqn:expressiong}
g_{\alpha,\beta}(s,t) = \frac{L_{\infty}\left(f,\frac{1}{2}+\alpha+s+it\right)L_{\infty}\left(f,\frac{1}{2}+\beta+s-it\right)}{L_{\infty}\left(f,\frac{1}{2}+\alpha+it\right)L_{\infty}\left(f,\frac{1}{2}+\beta-it\right)}, \hspace{1mm} V_{\alpha,\beta}(x,t) = \frac{1}{2i\pi}\int_{(1)}\frac{G(s)}{s}g_{\alpha,\beta}(s,t)x^{-s}ds 
\end{eqnarray*}\end{small}and
\begin{small}
\begin{eqnarray*} \label{eqn:expressionX}
X_{\alpha,\beta,t} = \frac{L_{\infty}\left(f,\frac{1}{2}-\alpha-it\right)L_{\infty}\left(f,\frac{1}{2}-\beta+it\right)}{L_{\infty}\left(f,\frac{1}{2}+\alpha+it\right)L_{\infty}\left(f,\frac{1}{2}+\beta-it\right)} .
\end{eqnarray*}
\end{small}
\end{lemma}
We do not write the proof of this lemma, which is essentially the same as theorem 5.3 of \cite{IK}. Nevertheless, it will be usefull to have good approximations of $X_{\alpha, \beta, t}$, $g_{\alpha,\beta}(s,t)$ and $V_{\alpha,\beta}(x,t)$.
\begin{lemma}\label{lem:estimationsdiverses} For large $t$ and for  $s\ll t^{\varepsilon}$ in any vertical strip, we have
\begin{eqnarray}
X_{\alpha, \beta, t} & = & \left(\frac{t\sqrt{N}}{2\pi}  \right)^{-2(\alpha+\beta)}\left( 1+\frac{i(\alpha^2-\beta^2)}{t}+O\left(  \frac{1}{t^2}  \right) \right)  \label{eqn:stirling} \\
\mbox{ and }g_{\alpha,\beta}(s,t) & = & \left( \frac{t\sqrt{N}}{2\pi} \right)^{2s}\left(1+O\left(\frac{|s^2|}{t}\right) \right) \label{eqn:stirling2}  .  
\end{eqnarray}
In addition, for each integer $j\geq 0$ and for all real number $A>0$, we have 
\begin{eqnarray} t^j\frac{\partial^j}{\partial t^j}V_{\alpha,\beta}(x,t) \ll_{A, j}\left( 1+ \frac{|x|}{t^2} \right)^{-A} . \label{eqn:estiV} \end{eqnarray}
\end{lemma}
Proof: We may write 
\begin{small}\begin{eqnarray*} X_{\alpha, \beta, t} & = & \left(\frac{\sqrt{N}}{2\pi}  \right)^{-2(\alpha+\beta)} \frac{\Gamma\left( \frac{k}{2}-\alpha-it \right)\Gamma\left( \frac{k}{2}-\beta+it \right)}{\Gamma\left( \frac{k}{2}+\alpha+it \right)\Gamma\left( \frac{k}{2}+\beta-it \right)} \\
 \mbox{ and }g_{\alpha,\beta}(s,t)& = & \left( \frac{\sqrt{N}}{2\pi} \right)^{2s}\frac{\Gamma\left( \frac{k}{2}+\alpha+s+it \right) \Gamma\left( \frac{k}{2}+\beta+s-it \right)}{\Gamma\left( \frac{k}{2}+\alpha+it \right) \Gamma\left( \frac{k}{2}+\beta-it \right)} .
 \end{eqnarray*} \end{small}Then, the first part of this lemma comes from the following Stirling formula with $s=\sigma+i\tau$ in any vertical strip  
\begin{small}$$ \Gamma(s)=\sqrt{2\pi}|\tau|^{\sigma-\frac12}e^{-\frac{\pi}{2}|\tau|}e^{i\left(\tau \ln |\tau| -\tau +\frac{\pi}{2}(\sigma-1/2)sgn(\tau) \right)}\left(1-i\frac{\left(\sigma-\frac12\right)^2-\frac{1}{12}}{2\tau} +O\left(\frac{1}{\tau^2}\right)\right).$$ \end{small}We refer to (\cite{Te}, corollaire 0.13). To prove (\ref{eqn:estiV}),  we move the integration line far to the right on $\Re s= A$ and by (\ref{eqn:stirling2}), we obtain the desired bound if $t^2\ll x$. In the case $x\ll t^2$, the result comes easily by trivial bounds.
\begin{flushright}
$\Box$
\end{flushright}
Thanks to the previous fonctional equation, we may split $I_f(\alpha,\beta)$ as a sum of diagonal terms and off-diagonal terms. Precisely, opening the mollifier $\psi$,  we may write
\begin{eqnarray} 
I_f(\alpha,\beta)& = &  \sum_{a,b\leq M } \frac{\mu_f(a)\mu_f(b)}{\sqrt{ab}}P\left(\frac{\ln M/a}{\ln M}  \right)P\left(\frac{\ln M/b}{\ln M}  \right)\left[I_{a,b}^{D_1}(\alpha,\beta)+  I_{a,b}^{D_2}(\alpha,\beta)\right. \label{eqn:expressionmomentdordre2twisted} \\
& & \hspace{8cm} \left. + I_{a,b}^{ND_1}(\alpha,\beta) + I_{a,b}^{ND_2}(\alpha,\beta)\right]\nonumber 
\end{eqnarray}
with
\begin{eqnarray*}
I_{a,b}^{D_1}(\alpha,\beta) & = & \sum_{am=bn}\frac{\lambda_f(m)\lambda_f(n)}{m^{\frac12+\alpha}n^{\frac12+\beta}}\int_{-\infty}^{+\infty}w(t)V_{\alpha, \beta}(mn,t)dt, \\
I_{a,b}^{D_2}(\alpha,\beta)& =&  \sum_{am=bn}\frac{\lambda_f(m)\lambda_f(n)}{m^{\frac12-\beta}n^{\frac12-\alpha}}\int_{-\infty}^{+\infty}w(t)X_{\alpha,\beta,t}V_{-\beta,-\alpha}(mn,t)dt, \\
I_{a,b}^{ND_1}(\alpha,\beta)& = & \sum_{am \neq bn}\frac{\lambda_f(m)\lambda_f(n)}{m^{\frac12+\alpha}n^{\frac12+\beta}}\int_{-\infty}^{+\infty}w(t)\left(\frac{am}{bn}\right)^{-it}V_{\alpha, \beta}(mn,t)dt, \\
I_{a,b}^{ND_2}(\alpha,\beta)& =&  \sum_{am\neq bn}\frac{\lambda_f(m)\lambda_f(n)}{m^{\frac12-\beta}n^{\frac12-\alpha}}\int_{-\infty}^{+\infty}w(t)X_{\alpha,\beta,t}\left(\frac{am}{bn}\right)^{-it}V_{-\beta,-\alpha}(mn,t)dt .
\end{eqnarray*}

\subsection{Evaluation of the off-diagonal term} \label{subsec:ODT}

In this part, we evaluate the size of the off-diagonal term. Precisely, we prove the following proposition. 
\begin{proposition} \label{pro:estimationtermeoffdiagonal} If $0<\nu< \frac{1-2\theta}{4+2\theta}$ and if  $\alpha, \beta$ are complex numbers satisfying  $\alpha, \beta \ll L^{-1}$ then there exists $\varepsilon >0$ such that 
\begin{small} 
$$\sum_{a,b\leq M } \frac{\mu_f(a)\mu_f(b)}{\sqrt{ab}}P\left(\frac{\ln M/a}{\ln M}  \right)P\left(\frac{\ln M/b}{\ln M}  \right)\left[ I_{a,b}^{ND_1}(\alpha,\beta) + I_{a,b}^{ND_2}(\alpha,\beta)\right] \ll  T^{1-\varepsilon}.$$
\end{small}
\end{proposition}
The main tool of the proof of this proposition is a theorem  about shifted convolution sums on average.

\subsubsection{Initial lemmas}
In order to prove the previous proposition, we begin to get rid of some harmless terms occurring in the definition of \begin{small}$I_{a,b}^{ND_1}(\alpha,\beta) $\end{small}.
\begin{lemma} Let $\varepsilon>0$, $0<\gamma<1$,  $\alpha, \beta \ll L^{-1}$ be complex numbers and  $a, b \leq T^\nu$ be positive integers. Then, for all real number $A>0$, we have
\begin{small}
\begin{eqnarray} \label{eqn:TND1}
I_{a,b}^{ND_1}(\alpha,\beta) =\sum_{\substack{am \neq bn \\ mn\ll T^{2+\varepsilon} \\ \left|\frac{am}{bn}-1 \right|\ll T^{-\gamma}}}\frac{\lambda_f(m)\lambda_f(n)}{m^{\frac12+\alpha}n^{\frac12+\beta}}\int_{-\infty}^{+\infty}w(t)\left(\frac{am}{bn}\right)^{-it}V_{\alpha, \beta}(mn,t)dt +O\left(T^{-A} \right). 
\end{eqnarray}
\end{small}
\end{lemma}
Proof:
Firstly, by (\ref{eqn:estiV}) with $j=0$, we get for all real number $A>0$
\begin{small} $$ \int_{-\infty}^{+\infty}w(t)\left(\frac{am}{bn}\right)^{-it}V_{\alpha, \beta}(mn,t)dt \ll_A T\left(\frac{T^2}{mn}\right)^A . $$ \end{small}As a consequence, for $A>1/2$, since \begin{small}$|\lambda_f(n)|\leq \tau(n) \ll n^\epsilon$\end{small} and $\alpha, \beta \ll L^{-1}$, we may write
\begin{small}
\begin{align*}
\sum_{\substack{am \neq bn \\ mn > T^{2+\varepsilon}}} \frac{\lambda_f(m)\lambda_f(n)}{m^{\frac12+\alpha}n^{\frac12+\beta}}\int_{-\infty}^{+\infty}w(t)\left(\frac{am}{bn}\right)^{-it}V_{\alpha, \beta}(mn,t)dt  & \ll T^{1+2A} \sum_{mn > T^{2+\varepsilon}}\frac{\tau(m)\tau(n)}{m^{\frac12+A+\Re \alpha}n^{\frac12+A+\Re \beta}} \\ 
& \ll T^{1+2A} \sum_{h > T^{2+\varepsilon}} \frac{1}{h^{\frac12 +A -\epsilon }}\\ & \ll T^{2-A\varepsilon} \\ & \ll T^{-A} .
\end{align*}
\end{small}Then, using (\ref{eqn:estiV}) for each integer $j$ and since $w^{(j)}(t) \ll \Delta^{-j}$, for all real number $A>0$, we have, uniformly with respect to $x$, the following bound
\begin{small} 
$$ \frac{\partial^j}{\partial t^j}\left[w(t)V_{\alpha,\beta}(x,t)\right] \ll \Delta^{-j}\left(\frac{T^2}{|x|} \right)^{A} .  $$ 
\end{small}Hence, if $am\neq bn$, thanks to $j$ integrations by parts, we get
\begin{small}
\begin{eqnarray*} \int_{-\infty}^{+\infty}w(t)\left(\frac{am}{bn}\right)^{-it}V_{\alpha, \beta}(mn,t)dt & = & \frac{1}{\left(i\ln\frac{am}{bn}\right)^j}\int_{-\infty}^{+\infty}\left(\frac{am}{bn}\right)^{-it}\frac{\partial^j}{\partial t^j}\left[w(t)V_{\alpha,\beta}(x,t)\right] dt \\
  & \ll &\frac{T}{\Delta^j \left| \ln\frac{am}{bn} \right|^j}\left(\frac{T^2}{mn} \right)^{A} .  \end{eqnarray*}
\end{small}Therefore, with $A=1/2+\max{\left\{\Re \alpha,\Re \beta\right\}} + \delta$ and $\delta>0$, using the lower bound $x/2 \leq \ln(1+x)$ for $0<x<1$, we  have
\begin{small}
$$
\sum_{\substack{ am \neq bn \\ \left|\frac{am}{bn}-1 \right|> T^{-\gamma} }} \frac{\lambda_f(m)\lambda_f(n)}{m^{\frac12+\alpha}n^{\frac12+\beta}}\int_{-\infty}^{+\infty}w(t)\left(\frac{am}{bn}\right)^{-it}V_{\alpha, \beta}(mn,t)dt \ll \frac{T^{1+2A+j\gamma}}{\Delta^j}\sum_{m,n\geq 1}\frac{\tau(m)\tau(n)}{m^{\frac12+\Re \alpha+A}n^{\frac12+\Re \beta+A}}.
$$
\end{small}Since $\gamma <1$, the result comes easily from a choice of large $j$.
\begin{flushright}
$\Box$
\end{flushright}
We introduce a dyadic partition of unity to the sums over $m$ and $n$. We fix an arbitrary smooth function $\rho : ]0,+\infty[ \rightarrow \R $,   which is compactly supported in $[1,2]$ and which satisfies
$$\sum_{\ell=-\infty}^{+\infty}\rho \left(2^{-\ell/2} x\right)=1 . $$
We may refer to \cite{Ha} (section 5) to build such a function. For each integer $\ell$, we define 
$$\rho_\ell(x)= \rho\left(\frac{x}{A_\ell} \right) \mbox{ with }A_\ell =2^{\ell/2}T^\gamma.$$
In order to study the asymptotic behaviour of $I_{a,b}^{ND_1}(\alpha,\beta)$, we consider the function $F_{h;\ell_1,\ell_2}$ which is defined by
\begin{eqnarray}
 F_{h;\ell_1,\ell_2}(x,y)= \frac{a^{\frac12+\alpha}b^{\frac12+\beta}}{x^{\frac12+\alpha}y^{\frac12+\beta}}\int_{-\infty}^{+\infty}w(t)\left(1+\frac{h}{y}\right)^{-it}V_{\alpha, \beta}\left(\frac{xy}{ab},t\right)dt\times \rho_{\ell_1}(x)\rho_{\ell_2}(y). 
 \end{eqnarray}
 \begin{lemma} \label{lem:simpl} Let $\varepsilon>0$, $0<\gamma<1$,  $\alpha, \beta \ll L^{-1}$ be complex numbers and  $a, b \leq T^\nu$ be positive integers. Then, for all real number $A>0$, we have
\begin{small}
\begin{eqnarray} 
I_{a,b}^{ND_1}(\alpha,\beta) = \sum_{\substack{ A_{\ell_1}A_{\ell_2} \ll abT^{2+\varepsilon} \\ A_{\ell_1} \asymp A_{\ell_2} \\ A_{\ell_1}, A_{\ell_2} \gg T^{\gamma} }} \sum_{0<|h|\ll T^{-\gamma}\sqrt{A_{\ell_1}A_{\ell_2}}} \sum_{am-bn=h} \lambda_f(m)\lambda_f(n) F_{h;\ell_1,\ell_2}(am,bn)
 +O\left(T^{-A} \right). 
\end{eqnarray}
\end{small}
\end{lemma}
Proof: For more convenience, we define
\begin{small}
$$ H(x,y)= \frac{a^{\frac12+\alpha}b^{\frac12+\beta}}{x^{\frac12+\alpha}y^{\frac12+\beta}}\int_{-\infty}^{+\infty}w(t)\left(\frac{x}{y}\right)^{-it}V_{\alpha, \beta}\left(\frac{xy}{ab},t\right)dt. $$
\end{small}From the previous lemma and using the partition of  unity, we may write
\begin{small}
$$I_{a,b}^{ND_1}(\alpha,\beta) = \sum_{\ell_1,\ell_2} \sum_{h\neq 0} \sum_{\substack{am - bn = h \\ mn\ll T^{2+\varepsilon} \\ \left|\frac{am}{bn}-1 \right|\ll T^{-\gamma}}}\lambda_f(m)\lambda_f(n)H(am,bn)\rho_{\ell_1}(am)\rho_{\ell_2}(bn) +O\left(T^{-A} \right). $$
\end{small}First, if $|h| \geq \sqrt{A_{\ell_1}A_{\ell_2}}T^{-\gamma}$ then 
\begin{small}
$$\max\left\{\left|\frac{am}{bn}-1\right|,\left|\frac{bn}{am}-1\right|\right\}^{2} \geq \left|\frac{am}{bn}-1\right|\left|\frac{bn}{am}-1\right|=\frac{h^2}{ambn} \asymp  \frac{h^2}{A_{\ell_1}A_{\ell_2}} \geq T^{-2\gamma} . $$
\end{small}Secondly, if $|\ell_1-\ell_2|\geq 3$, for instance if $\ell_1-\ell_2\geq 3$ then 
\begin{small}
$$ \frac{am}{bn}-1 \geq \frac{2^{\frac{\ell_1-\ell_2}{2}}}{2}-1 \geq \sqrt{2}-1 \gg 1 .$$
\end{small}Thus we may assume $A_{\ell_1}\asymp A_{\ell_2}$.
Thirdly, if \begin{small} $A_{\ell_2} \leq T^{\gamma} $  \end{small} 
then
\begin{small}
$$\left|\frac{am}{bn}-1  \right|=\frac{h}{bn} \geq \frac{h}{2A_{\ell_2}}\gg T^{-\gamma}.$$
\end{small}Therefore, we may assume \begin{small} $A_{\ell_2}\geq T^\gamma $  \end{small} and, in the same way, \begin{small}$A_{\ell_1}\geq T^\gamma $ \end{small}. Finally, since $am-bn=h$, we get $H(am,bn)\rho_{\ell_1}(am)\rho_{\ell_2}(bn)=F_{h;\ell_1,\ell_2}(am,bn)$.
\begin{flushright}
$\Box$
\end{flushright}

\subsubsection{Shifted convolution sums}
The core of the proof of our theorem is the following bound, which is a generalisation of theorem 2 of \cite{Bl2}, about shifted convolution sums on average. We define $\theta$ as the exponent in the Ramanujuan-Petersson conjecture, which claims
\begin{eqnarray}\label{eqn:defRP}
| \lambda(n) | \leq \tau(n)n^\theta
\end{eqnarray}
for eigenvalues $\lambda(n)$ of the Hecke operator $T_n$ acting on the space of weight $0$ Maa$\mbox{\ss}$ cusp forms of level $N$.
\begin{theorem} \label{th:thBlomergenenmoyenne}Let $\ell_1$, $\ell_2$, $H$ and $h_1$ be positive integers. Let $M_1$, $M_2$, $P_1$, $P_2$ be real numbers greater than $1$. Let $\{ g_h \}$ be a family of smooth functions, supported on $[M_1,2M_1]\times [M_2,2M_2]$ such that $||g_h^{(ij)}||_{\infty} \ll_{i,j}(P_1/M_1)^i(P_2/M_2)^j$ for all $i,j \geq 0$. Let $(a(h))$ be a sequence of complex numbers such that
$$ a(h) \neq 0 \Rightarrow h \leq H, \; h_1 |h \mbox{ et } \left( h_1,\frac{h}{h_1}\right) =1. $$
If  $\ell_1M_1 \asymp \ell_2M_2 \asymp A$ and if there exists $\epsilon >0$ such that \begin{eqnarray}H \ll  \frac{A}{\max\{P_1,P_2\}^2}\frac{1}{(\ell_1\ell_2M_1M_2P_1P_2)^{\epsilon}}, 
\end{eqnarray}
then, for all real number $\varepsilon >0$, we have
\begin{small}
\begin{multline*}
\sum_{h=1}^H a(h)\sum_{\substack{m_1,m_2 \geq 1 \\ \ell_1m_1-\ell_2m_2 = h}} \lambda_f(m_1)\overline{\lambda_f(m_2)}g_h(m_1,m_2) \\
\ll A^{\frac12}h_1^\theta \|a\|_2(P_1+P_2)^{\frac32} \left[\sqrt{P_1+P_2}+\left( \frac{A}{\max\{ P_1,P_2\}}\right)^\theta \left( 1+ \sqrt{\frac{(h_1,\ell_1\ell_2)H}{h_1\ell_1\ell_2}}\right) \right](\ell_1\ell_2M_1M_2P_1P_2H)^\varepsilon .
\end{multline*}
\end{small}
\end{theorem}
Proof: This proof is a direct generalisation of the proof of theorem 2 in \cite{Bl2}, that is why we only give an outline of the proof of Blomer. In order to simplify, we set
\begin{small}
$$
\Sigma(\ell_1,\ell_2,H,a)  = \sum_{h=1}^H a(h)\sum_{\substack{m_1,m_2 \geq 1 \\ \ell_1m_1-\ell_2m_2 = h}} \lambda_f(m_1)\overline{\lambda_f(m_2)}g_h(m_1,m_2) .
$$\end{small}We set  also $\Delta = \min \left\{\frac{P_1}{\ell_1M_1},\frac{P_2}{\ell_2M_2}\right\}$,  let $Q \geq 1/\Delta$ be a very large parameter and set $\delta =1/Q$ such that $\delta \leq \Delta$. Let $\phi$ be a smooth function compactly supported in $[-\Delta^{-1},\Delta^{-1}]$, with $\phi(0)=1$ and such that, for all integer $j$, $\| \phi^{(j)} \|_{\infty} \ll_j \Delta^j$. Then, we introduce
   $$W_h(x,y)=g_h(x,y)\phi(\ell_1x-\ell_2y-h) .$$
In addition, we introduce another smooth function   $w :\R \rightarrow \R$ compactly supported in  $[Q,2Q]$ and such that $\| w^{(j)} \|_{\infty} \ll_j Q^{-j}$. If $\varphi$ denotes Euler's phi function, let \begin{eqnarray*}\Lambda = \sum_{q\equiv 0[N\ell_1\ell_2]}w(q)\varphi(q) \asymp \frac{Q^2}{N\ell_1\ell_2}. \end{eqnarray*}
As a result, we can rewrite\begin{small}
\begin{align*}\Sigma(\ell_1,\ell_2,H,a) = \sum_{h=1}^H a(h)\int_{0}^{1}\sum_{m_1,m_2 \geq 1} \l(m_1)\overline{\l(m_2)}e(\ell_1m_1\alpha)e(-\ell_2m_2\alpha)W_h(m_1,m_2)e(-h\alpha)d\alpha. 
\end{align*}\end{small}By the Jutila's circle method (see \cite{Bl2}, lemma 3.1), we build an approximation $\tilde{I}$ to the characteristic function on $[0,1]$, which splits 
 the $\alpha$-integral in two parts, according to whether $\alpha$ is in a minor arc or a major arc. We easily estimate the contribution of minor arcs using our bound for the $L^2([0,1])$-norm of $1-\tilde{I}$ and we can write the contribution of major arcs of the shape
\begin{multline*} 
 \frac{1}{2\delta\Lambda}\sum_{h=1}^H a(h)\sum_{q\equiv 0 [N\ell_1\ell_2]}w(q){\sum_{d(q)}}^* \int_{-\delta}^{\delta} \sum_{m_1,m_2 \geq 1 } \l(m_1)\overline{\l(m_2)} \nonumber \\  e\left(\ell_1m_1\left(\frac{d}{q}+\eta\right) \right)e\left(-\ell_2m_2\left(\frac{d}{q}+\eta\right) \right)W_h(m_1,m_2)e\left(-h\left(\frac{d}{q}+\eta\right) \right)d\eta .
\end{multline*} 
We transform short sums of exponentials in long sums of exponentials by means of a Voronoï summation formula (see \cite{Bl2}, lemma 2.2 or \cite{KMV}, appendix A). If $S(m,n;q)$ denotes the classical Kloosterman sum, the major arcs contribution becomes 
\begin{multline*}
   \frac{1}{2\delta\Lambda}\sum_{h=1}^H a(h)\sum_{q\equiv 0 [N\ell_1\ell_2]}w(q) \int_{-\delta}^{\delta} e(-\eta h)\\  \times \sum_{m_1,m_2 \geq 1 } \l(m_1)\overline{\l(m_2)} S(-h,\ell_2m_2-\ell_1m_1;q)G_{q,h,\eta}(m_1,m_2) d\eta
\end{multline*}where
\begin{multline*}
 G_{q,h,\eta}(x_1,x_2)=\frac{4\pi^2\ell_1\ell_2}{q^2}\int_{0}^{+\infty}\int_{0}^{+\infty}W_h(t_1,t_2)e(\ell_1t_1\eta-\ell_2t_2\eta)\\ \times J_{k-1}\left( \frac{4\pi \ell_1 \sqrt{x_1t_1}}{q}\right)J_{k-1}\left( \frac{4\pi \ell_2 \sqrt{x_2t_2}}{q}\right)dt_1dt_2.\end{multline*}
We split this sum according to the value of $\ell_2m_2-\ell_1m_1$. It comes a diagonal contribution when $\ell_1m_1=\ell_2m_2$, which is easily bounded, and an off-diagonal contribution when $\ell_2m_2\neq \ell_1m_1$. It remains to estimate this off-diagonal contribution, which can be rewritten
\begin{multline*}
 \frac{\pi \ell_1\ell_2}{2\Lambda}\int_{-\delta}^{\delta}\int_0^{+\infty}\int_0^{+\infty}\sum_{h=1}^H a_{\eta,t_1,t_2}(h)\sum_{r\neq 0}\sum_{\substack{1\leq m_1 \leq \mathcal{M}_1 \\ 1\leq m_2 \leq \mathcal{M}_2 \\ \ell_2m_2-\ell_1m_1=r} }\l(m_1)\overline{\l(m_2)} \\ 
 \times \sum_{q\equiv 0 [N\ell_1\ell_2]}\frac{S(-h,r;q)}{q} \Phi_{t_1,t_2}\left( \frac{4\pi \sqrt{h|r|}}{q};m_1,m_2,r,h\right)dt_1dt_2 d\eta .
\end{multline*}  
 where
 \begin{align*}
  a_{\eta,t_1,t_2}(h) &  = a(h)e(-\eta h )W_h(t_1,t_2)e(\ell_1t_1\eta -\ell_2t_2\eta) \\
\mbox{and } \Phi_{t_1,t_2}(x;m_1,m_2,r,h) & =  \frac{Qx}{\sqrt{h|r|}}J_{k-1}\left(\frac{x\ell_1 \sqrt{m_1t_1}}{\sqrt{h|r|}} \right)J_{k-1}\left(\frac{x\ell_2 \sqrt{m_2t_2}}{\sqrt{h|r|}} \right)w\left(\frac{4\pi\sqrt{h|r|}}{x}\right)w_1(h)
 \end{align*}
and with $$\mathcal{M}_1 =\frac{Q^2P_1^2k^2}{\ell_1^2M_1}(\ell_1\ell_2M_1M_2P_1P_2)^\varepsilon \mbox{ and }\mathcal{M}_2 =\frac{Q^2P_2^2k^2}{\ell_2^2M_2}(\ell_1\ell_2M_1M_2P_1P_2)^\varepsilon.$$
The asymptotic behaviour of $G_{q,h,\eta}$ allows us to restrict sums over  $m_1$ and $m_2$ to sums over $m_1\leq \mathcal{M}_1$ and $m_2\leq \mathcal{M}_2$. Applying the Kuznetsov trace formula (see \cite{Bl2}, lemma 2.4 or \cite{DI}, theorem 1), we decompose this off-diagonal term as a sum of three terms: the contribution of discrete spectrum, the contribution of continuous spectrum and the contribution of holomorphic cusp forms. All of them will be evaluated by means of large sieve inequalities (see \cite{Bl2}, lemma 2.5 or \cite{DI}, theorem 2). In the discrete spectrum, it may have exceptional eigenvalues even thought the Selberg's conjecture predicts that such exceptional eigenvalues do not exist. We prove that the contribution of non-exceptional eigenvalues, called real spectrum, is bounded by 
$$A^{1/2}\|a\|_2h_1^\theta (P_1+P_2)^{\frac32}\left( \sqrt{P_1+P_2}+\sqrt{\frac{H(h_1,\ell_1\ell_2)}{h_1\ell_1\ell_2}} \right)Q^\varepsilon$$ whereas we bound the contribution of exceptional eigenvalues by 
$$A^{\frac{1}{2}+\theta}h_1^\theta \|a\|_2 (P_1+P_2)^{\frac32 - \theta} \left(1+\sqrt{\frac{H(h_1,\ell_1\ell_2)}{h_1\ell_1\ell_2}} \right)Q^\varepsilon.$$
In addition, we show that the contribution of the continuous spectrum and of holomorphic cusp forms are bounded by 
$$A^{1/2}\|a\|_2 (P_1+P_2)^{3/2}\left( \sqrt{P_1+P_2}+\sqrt{\frac{H(h_1,\ell_1\ell_2)}{h_1\ell_1\ell_2}} \right)Q^\varepsilon,$$
which is smaller than the estimate of the real spectrum. The result comes easily from these last three estimates.
\begin{flushright}
$\Box$
\end{flushright}
\begin{remark} In \cite{Bl2}, a factor $\left(\frac{LM}{HP}\right)^\theta$ appears, which could become $\left(\frac{A}{H(P_1+P_2)}\right)^\theta$ in our theorem, and which comes from the contribution of possibly exceptional eigenvalues in the discrete spectrum. However, we only find $\left(\frac{A}{P_1+P_2}\right)^\theta$. 
\end{remark}
\subparagraph*{}
Let us determine the required bounds for the test function in our case.
\begin{lemma} \label{lem:deriveepartielle} Let  $\alpha, \beta \ll L^{-1}$ be complex numbers and let $\sigma$ be any positive real number. For all non negative integers $i$ and $j$, we have
\begin{small}
$$
x^{i}y^j \frac{\partial^{i+j}F_{h;\ell_1,\ell_2}}{\partial x^{i}\partial y^j}(x,y) \ll_{i,j}  \left(\frac{a}{A_{\ell_1}}\right)^{\frac12+\Re \alpha+\sigma}\left(\frac{b}{A_{\ell_2}}\right)^{\frac12+\Re \beta+\sigma} T^{1+2\sigma} (\ln T)^j
$$
and the implicit constant does not depend on $h$.
\end{small}
\end{lemma}
Proof: Let $P_0(Y)=1$ and $P_j(Y)=\prod_{\ell=0}^{j-1}(Y-\ell)$ for $j\geq 1$. Let  $\Re s=\sigma >0$, then we have
\begin{small}
\begin{eqnarray*}
y^j\frac{\partial^j}{\partial y^j}\left(\int_{\R}w(t)\left(1+\frac{h}{y}\right)^{-it}g_{\alpha,\beta}(s,t)dt \right)  & = & \int_\R w(t)g_{\alpha,\beta}(s,t)\left(1+\frac{h}{y}\right)^{-it}Q_j(t)dt \\
 & = & \int_\R   \frac{ \left(1+\frac{h}{y}\right)^{-it}}{\left[i \ln \left(1+\frac{h}{y} \right) \right]^j}          \frac{\partial^j}{\partial t^j}\left[ w(t)g_{\alpha,\beta}(s,t)Q_j(t)\right]dt
\end{eqnarray*}
\end{small}where \begin{small} $Q_j(t)=\sum_{r=0}^j \binom{j}{r}P_r(it)P_{j-r}(-it)\left(1+\frac{h}{y}\right)^{-(j-r)}$ \end{small}.
Since $Q_j^{(r)}(t) \ll \left|\frac{h/y}{1+h/y}\right|^j t^{j-r}$ for $r\leq j$ and $\frac{\partial^r}{\partial t^r}\left[ w(t)g_{\alpha,\beta}(s,t)\right]dt \ll T^{2\sigma}\Delta^{-r}$, we get 
\begin{small}$$ \frac{\partial^j}{\partial t^j}\left[ w(t)g_{\alpha,\beta}(s,t)Q_j(t)\right] \ll_{s,j} T^{2\sigma}\left|\frac{h/y}{1+h/y}\right|^j (\ln T)^j. $$\end{small}and the implicit constant depends polynomially on $s$. Since $h/y \ll T^{-\gamma}$ in our range of summation, then  
\begin{small}
\begin{eqnarray} \label{eqn:boundderpar}
y^j\frac{\partial^j}{\partial y^j}\left(\int_{\R}w(t)\left(1+\frac{h}{y}\right)^{-it}g_{\alpha,\beta}(s,t)dt \right) \ll_{s,j} T^{1+2\sigma} (\ln T)^j \end{eqnarray}
\end{small}and the implicit constant does not depend on $h$. Writing 
\begin{small}
$$F_{h;\ell_1,\ell_2}(x,y) = \frac{\rho_{\ell_1}(x)\rho_{\ell_2}(y)}{2i\pi}\int_{(\sigma)}\frac{G(s)}{s}\left(\frac{a}{x} \right)^{\frac12+\alpha+s}\left(\frac{b}{y} \right)^{\frac12+\beta+s}\int_{\R}w(t)\left(1+\frac{h}{y}\right)^{-it}g_{\alpha,\beta}(s,t)dt ds,  $$
\end{small}the result comes easily from the Leibniz formula and the bound (\ref{eqn:boundderpar}).
\begin{flushright}
$\Box$
\end{flushright}
\begin{remark} The trivial bound for shifted convolution sums, namely taking absolute values and applying the Ramanujuan-Petersson bound on average, is given, for all $\varepsilon>0$, by 
\begin{align}
\sum_{\ell_1 m_1 -\ell_2 m_2= h} \lambda_f(m_1)\lambda_f(m_2) g_h(m_1,m_2) \ll_\varepsilon \min \{M_1,M_2\} (M_1M_2)^\varepsilon . \label{eqn:trivialbound}
\end{align}
The trivial bound (\ref{eqn:trivialbound}) and lemma \ref{lem:deriveepartielle} imply the following corollary.
\end{remark}
\begin{corollary} For all $\varepsilon >0$, we have
\begin{small}$$I_{a,b}^{ND_1}(\alpha,\beta) \ll_\varepsilon \min \{a,b\}T^{1+\varepsilon}. $$\end{small}
\end{corollary}
\begin{remark} This trivial bound of \begin{small}$I_{a,b}^{ND_1}(\alpha,\beta)$\end{small} fails to prove Proposition \ref{pro:estimationtermeoffdiagonal}. In other words, taking care of the oscillations of the Hecke eigenvalues is required to proving Proposition \ref{pro:estimationtermeoffdiagonal}. In a first time, we apply the following bound, proved by Blomer in \cite{Bl}, about shifted convolution sums.
\end{remark}
\begin{theorem}[Blomer, \cite{Bl}] \label{th:thBlomer}
 Let $\varepsilon >0$, and let $\ell_1$, $\ell_2$, $h$ be positive integers. Let $M_1$, $M_2$, $P_1$, $P_2$ be real numbers greater than $1$. Let $g_h$ be a smooth function, supported on $[M_1,2M_1]\times [M_2,2M_2]$ such that $||g_h^{(ij)}||_{\infty} \ll_{i,j}(P_1/M_1)^i(P_2/M_2)^j$ for all $i,j \geq 0$. Then
$$\sum_{\ell_1 m_1 -\ell_2 m_2= h} \lambda_f(m_1)\lambda_f(m_2) g_h(m_1,m_2) \ll_{\varepsilon,P_1,P_2,N,k} (\ell_1 M_1+\ell_2 M_2)^{1/2+\theta+\varepsilon}.  $$
This bound is uniform in $\ell_1$, $\ell_2$, $h$ and the dependence on $P_1$, $P_2$, $N$ and $k$ is polynomial.
\end{theorem}
\begin{remark}
Remembering the trivial bound (\ref{eqn:trivialbound}), this theorem agrees with the square-root cancellation philosophy.
\end{remark}
\subparagraph*{}
As a consequence, theorem \ref{th:thBlomer} and lemma \ref{lem:deriveepartielle} imply the following proposition, which gives a first admissible bound, and which will be improved in proposition \ref{pro:estiOD1bis}.  
\begin{proposition} \label{pro:estiOD1} Let $\alpha, \beta \ll L^{-1}$ be complex numbers and  $a, b$ be positive integers. For all $\varepsilon >0$, we have
\begin{small}
$$I_{a,b}^{ND_1}(\alpha,\beta) \ll_{\varepsilon} (ab)^{\frac34+\frac{\theta}{2}}T^{\frac12+\theta+\varepsilon} .$$
\end{small}
\end{proposition}
Proof: By theorem \ref{th:thBlomer}, lemma \ref{lem:deriveepartielle} gives 
\begin{small}
$$ \sum_{am-bn=h} \lambda_f(m)\lambda_f(n) F_{h;\ell_1,\ell_2}(am,bn) \ll  \left(\frac{a}{A_{\ell_1}}\right)^{\frac12+\Re \alpha+\sigma}\left(\frac{b}{A_{\ell_2}}\right)^{\frac12+\Re\beta+\sigma} T^{1+2\sigma} (\ln T)^\kappa (A_{\ell_1}+A_{\ell_2})^{\frac12+\theta+\varepsilon} $$
\end{small}where $\kappa$ is a constant. Thus, thanks to lemma \ref{lem:simpl} and with  $ \frac{1}{2}+\theta+\varepsilon-\Re \alpha-\Re \beta-2\sigma >0$, we get
\begin{small}
\begin{align*} I_{a,b}^{ND_1}(\alpha,\beta) & \ll   T^{1-\gamma+2\sigma}(\ln T)^\kappa \sum_{\substack{ A_{\ell_1}A_{\ell_2} \ll abT^{2+\varepsilon} \\ A_{\ell_1} \asymp A_{\ell_2} \\ A_{\ell_1}, A_{\ell_2} \gg T^{\gamma} }}\left(\frac{a}{A_{\ell_1}}\right)^{\frac12+\Re\alpha+\sigma}\left(\frac{b}{A_{\ell_2}}\right)^{\frac12+\Re\beta+\sigma} \sqrt{A_{\ell_1}A_{\ell_2}}(A_{\ell_1}+A_{\ell_2})^{1/2+\theta+\varepsilon} \\ 
& \ll     T^{1-\gamma+2\sigma}(ab)^{1/2+\sigma} (\ln T)^\kappa \sum_{T^\gamma \ll A_{\ell_1} \ll \sqrt{ab}T^{1+\varepsilon/2}} A_{\ell_1}^{\frac{1}{2}+\theta+\varepsilon-\Re \alpha-\Re \beta-2\sigma}   \\
& \ll  T^{1-\gamma }T^{\frac12+\theta+\varepsilon}(ab)^{\frac34+\frac{\theta}{2}}\sum_{T^\gamma \ll A_{\ell_1} \ll \sqrt{ab}T^{1+\varepsilon/2}}1.
\end{align*}
\end{small}Then, we obtain $\sum_{T^\gamma \ll A_{\ell_1} \ll \sqrt{ab}T^{1+\varepsilon/2}}1 = \sum_{1\leq 2^{\ell_1/2} \ll \sqrt{ab}T^{1-\gamma+\varepsilon/2}}1\ll \ln T$. Finally, the result easily comes from the choice $\gamma=1-\varepsilon$.
\begin{flushright}
$\Box$
\end{flushright}
\begin{remark} We are tempted to solve the shifted convolution problem on average (over $h$) and to take care of the resulting additional oscillations of the Hecke eigenvalues. For instance, using Theorem $6.3$ of \cite{Ric}, one can check that, for all $\varepsilon >0$, 
\begin{small}$$ I_{a,b}^{ND_1}(\alpha,\beta) \ll_\varepsilon (ab)^{\frac34+\frac{\theta}{2}}T^{\frac32+\theta+\varepsilon}.$$\end{small}It turns out that this bound is not admissible. It can be explained by the fact that the length of the $h$-sum is very small. That is why we need a bound for short sums of shifted convolution sums, which is given  by theorem \ref{th:thBlomergenenmoyenne}. 
\end{remark}
\begin{proposition}\label{pro:estiOD1bis} Let $\alpha, \beta \ll L^{-1}$ be complex numbers and  $a, b$ be positive integers. For all $\varepsilon >0$, we have
\begin{small}
$$I_{a,b}^{ND_1}(\alpha,\beta) \ll_{\varepsilon} (ab)^{\frac{1+\theta}{2}}T^{\frac12+\theta+\varepsilon} .$$
\end{small}
\end{proposition}
Proof: We apply theorem \ref{th:thBlomergenenmoyenne} with $H=T^{-\gamma}\sqrt{A_{\ell_1}A_{\ell_2}}$, $h_1=1$ and $a(h)=\begin{cases}1 \mbox{ if }h\leq H, \\ 0 \mbox{ otherwise}. \end{cases} $
Therefore, from lemma \ref{lem:simpl}, we get
\begin{small}
\begin{align*}
& I_{a,b}^{ND_1}(\alpha,\beta)\\  \ll &   \sum_{\substack{ A_{\ell_1}A_{\ell_2} \ll abT^{2+\varepsilon} \\ A_{\ell_1} \asymp A_{\ell_2} \\ A_{\ell_1}, A_{\ell_2} \gg T^{\gamma} }} \left(\frac{a}{A_{\ell_1}}\right)^{\frac12+\Re \alpha+\sigma}\left(\frac{b}{A_{\ell_2}}\right)^{\frac12+\Re \beta+\sigma} T^{1+2\sigma+\varepsilon}\sqrt{A_{\ell_1}H} \left[\sqrt{\ln T}+ \left(\frac{A_{\ell_1}}{\ln T} \right)^\theta \left( 1+ \sqrt{\frac{H}{ab}} \right) \right] \\
 \ll & (ab)^{\frac12+\theta}T^{1+2\sigma +\varepsilon -\frac{\gamma}{2}}\sum_{T^\gamma \ll A_\ell \ll \sqrt{ab}T^{1+\varepsilon}}A_\ell^{\theta -(\Re \alpha + \Re \beta +2\sigma)} \ll (ab)^{\frac{1+\theta}{2}}T^{1-\frac{\gamma}{2}+\theta+\varepsilon}.
\end{align*}\end{small}Finally, the result comes from $\gamma= 1-\varepsilon$. 
\begin{flushright}
$\Box$
\end{flushright}
\begin{corollary} \label{pro:estiOD2} Let $\alpha, \beta \ll L^{-1}$ be complex numbers and  $a, b$ be positive integers. For all $\varepsilon >0$, we have
\begin{small}
$$I_{a,b}^{ND_2}(\alpha,\beta) \ll_\varepsilon (ab)^{\frac{1+\theta}{2}}T^{\frac12+\theta+\varepsilon} .$$
\end{small}
\end{corollary}
Proof: Set \begin{small}$$w_1(t)=w(t)\left(\frac{t\sqrt{N}}{2\pi}  \right)^{-2(\alpha+\beta)}\left( 1+\frac{i(\alpha^2-\beta^2)}{t}\right). $$\end{small}Thanks to relation (\ref{eqn:stirling}), we may write
\begin{small}
$$ I_{a,b}^{ND_2}(\alpha,\beta) = \sum_{am\neq bn}\frac{\lambda_f(m)\lambda_f(n)}{m^{\frac12-\beta}n^{\frac12-\alpha}}\int_{-\infty}^{+\infty}w_1(t)\left(\frac{am}{bn}\right)^{-it}V_{-\beta,-\alpha}(mn,t)dt +O\left(\frac{1}{T}\sum_{mn \ll T^{2+\varepsilon}} \frac{\left|\lambda_f(m)\lambda_f(n)\right|}{m^{\frac12-\Re \beta}n^{\frac12-\Re\alpha}}\right).$$
\end{small}The error term becomes $O(T^\varepsilon)$ and since $w_1$ satisfies (\ref{eqn:hypwa}), (\ref{eqn:hypwb}) and (\ref{eqn:hypwc}) we may apply  proposition \ref{pro:estiOD1bis} up to changing $w$ by $w_1$ and $(\alpha,\beta)$ by $(-\beta,-\alpha)$.
\begin{flushright}
$\Box$
\end{flushright}

\subsubsection{Proof of proposition \ref{pro:estimationtermeoffdiagonal}}
Using proposition \ref{pro:estiOD1bis} and corollary \ref{pro:estiOD2}, we trivially bound the off-diagonal term by 
\begin{small}
\begin{align*}
 \sum_{a,b\leq M } \frac{\mu_f(a)\mu_f(b)}{\sqrt{ab}}P\left(\frac{\ln M/a}{\ln M}  \right)P\left(\frac{\ln M/b}{\ln M}  \right)\left[ I_{a,b}^{ND_1}(\alpha,\beta) + I_{a,b}^{ND_2}(\alpha,\beta)\right] & \ll T^{\frac12+\theta+\varepsilon}  \sum_{a,b\leq M }(ab)^{\frac{\theta}{2}} \\ 
& \ll T^{\frac12+\theta+\varepsilon} T^{\nu\left(2+\theta\right)}.
 \end{align*}\end{small}
 
\noindent Thus, if $\nu< \frac{1-2\theta}{4+2\theta}$, the off-diagonal part of $I_f(\alpha,\beta)$ is bounded by $T^{1-\varepsilon}$. 

\subsection{Evaluation of the diagonal term} \label{subsec:EDT}

For $i=1$ or $i=2$, let  
\begin{eqnarray} \label{eqn:defIfD1}
I_{f}^{D_i}(\alpha, \beta) = \sum_{a,b\leq M } \frac{\mu_f(a)\mu_f(b)}{\sqrt{ab}}P\left(\frac{\ln M/a}{\ln M}  \right)P\left(\frac{\ln M/b}{\ln M}  \right)I_{a,b}^{D_i}(\alpha,\beta).
\end{eqnarray}
We consider the diagonal part $I_f^D(\alpha, \beta)$ of the mollified second moment. Thus, we have
\begin{eqnarray} \label{eqn:splitTD}
 I_{f}^{D}(\alpha, \beta) =  I_{f}^{D_1}(\alpha, \beta) +  I_{f}^{D_2}(\alpha, \beta).
 \end{eqnarray}
In this section, we prove the following proposition:
\begin{proposition} \label{pro:termmediag} Let $0<\nu<1$. For complex numbers $\alpha, \beta \ll L^{-1}$ such that $|\alpha+\beta | \gg L^{-1}$, we have 
$$I_{f}^{D}(\alpha, \beta ) = \widehat{w}(0)c(\alpha, \beta) +O(T(\ln L)^4/L)$$
where $ c(\alpha,\beta)$ is defined in relation \emph{(\ref{eqn:limiteM2})}.
\end{proposition}

\subsubsection{Initial lemmas}

\begin{lemma} \label{eqn:expressionproduiteulerienenO} Set $\Omega_{\alpha,\beta}$ the subset of vectors $(u,v,s)$ in $\C^3$ satisfying 
$$  \left\{ \begin{array}{l} \Re u+\Re v>-1/2, \\ \Re s>-\frac14-\frac{\Re\alpha}{2}-\frac{\Re\beta}{2},\\ \Re u+\Re s >-\frac12-\Re  \alpha,\\ \Re v+\Re s>-\frac12- \Re \beta . \end{array} \right. $$
 We have 
$$\sum_{\substack{a,b,m,n \geq 1 \\am=bn }}\frac{\mu_f(a)\mu_f(b)\lambda_f(m)\lambda_f(n)}{a^{\frac12+v}b^{\frac12+u}m^{\frac12+\alpha+s}n^{\frac12+\beta+s}} = \frac{L(f\times f, 1+\alpha+\beta+2s)L(f\times f, 1+u+v)}{L(f\times f, 1+\alpha+u+s)L(f\times f, 1+\beta+v+s)}A_{\alpha,\beta}(u,v,s)$$
where $A_{\alpha,\beta}(u,v,s)$ is given by an absolutely convergent Euler product on $\Omega_{\alpha,\beta}$.  
\end{lemma}
Proof : In order to simplify, set \begin{small}$$\mathcal{P}=\sum_{\substack{a,b,m,n \geq 1 \\am=bn }}\frac{\mu_f(a)\mu_f(b)\lambda_f(m)\lambda_f(n)}{a^{\frac12+v}b^{\frac12+u}m^{\frac12+\alpha+s}n^{\frac12+\beta+s}} .$$
\begin{normalsize}Using relation (\ref{eqn:multiplicativite}), for any prime number $p$ such that $p \nmid N$, we get,\end{normalsize}
$$
\sum_{\ell\geq 0}\frac{\lambda_f(p^\ell)\lambda_f(p^{\ell+1})}{p^{\ell s}} = \lambda_f(p)\sum_{\ell\geq 0} \frac{\lambda_f(p^\ell)^2}{p^{\ell s}}-\frac{1}{p^s}\sum_{\ell\geq 0}\frac{\lambda_f(p^{\ell+1})\lambda_f(p^\ell)}{p^{\ell s}} $$
\begin{normalsize}and\end{normalsize} $$ \sum_{\ell\geq 0}\frac{\lambda_f(p^\ell)\lambda_f(p^{\ell+2})}{p^{\ell s}} = \lambda_f(p^2)\sum_{\ell\geq 0} \frac{\lambda_f(p^\ell)^2}{p^{\ell s}}-\frac{\lambda_f(p)}{p^s}\sum_{\ell\geq 0}\frac{\lambda_f(p^{\ell+1})\lambda_f(p^\ell)}{p^{\ell s}} . $$
\begin{normalsize}Thus, we  deduce\end{normalsize}  
\begin{eqnarray}
 \sum_{\ell\geq 0}\frac{\lambda_f(p^\ell)\lambda_f(p^{\ell+1})}{p^{\ell s}}   & =  & \lambda_f(p) \left( 1+\frac{1}{p^s}\right)^{-1}\sum_{l\geq 0} \frac{\lambda_f(p^\ell)^2}{p^{\ell s}}  \label{eqn:relationmulti1} \\
 \sum_{\ell\geq 0}\frac{\lambda_f(p^\ell)\lambda_f(p^{\ell+2})}{p^{\ell s}}   & = & \left( 1+\frac{1}{p^s}\right)^{-1} \left(\lambda_f(p^2)-\frac{1}{p^s} \right)  \sum_{\ell\geq 0}\frac{\lambda_f(p^\ell)^2}{p^{\ell s}} . \label{eqn:relationmulti2} 
\end{eqnarray}
\begin{normalsize}In addition, since we write\end{normalsize}
\begin{eqnarray*}
\mathcal{P}= \prod_p \left(\sum_{\substack{\ell_1,\ell_2,\ell_3,\ell_4\geq 0\\ \ell_1+\ell_3=\ell_2+\ell_4}} \frac{\mu_f(p^{\ell_1})\mu_f(p^{\ell_2})\lambda_f(p^{\ell_3})\lambda_f(p^{\ell_4})}{p^{\ell_1(\frac12+v)}p^{\ell_2(\frac12+u)}p^{\ell_3(\frac12+\alpha+s)}p^{\ell_4(\frac12+\beta+s)}} \right) 
\end{eqnarray*}
\begin{normalsize}using  relations (\ref{eqn:proprietemu1}) and  (\ref{eqn:proprietemu2}), we get\end{normalsize}
\begin{eqnarray*}
 & & \mathcal{P} = \prod_{p\nmid N} \left[  \left(1+\frac{\lambda_f(p)^2}{p^{1+u+v}}+\frac{1}{p^{2(1+u+v)}}\right)\sum_{\ell\geq 0} \frac{\lambda_f(p^\ell)^2}{p^{\ell(1+\alpha+\beta+2s)}} \right. \\   
 & & \hspace{47mm} -\lambda_f(p)\left(\frac{1}{p^{1+v+\beta+s}} 
 + \frac{1}{p^{1+u+\alpha+s}}\right)\left(1+\frac{1}{p^{1+u+v}}\right)\sum_{\ell\geq 0}\frac{\lambda_f(p^\ell)\lambda_f(p^{\ell+1})}{p^{\ell(1+\alpha+\beta+2s)}} \\
 & & \hspace{68mm}\left.
 +\left(\frac{1}{p^{2(1+v+\beta+s)}}+\frac{1}{p^{2(1+u+\alpha+s)}}\right)\sum_{\ell\geq 0}\frac{\lambda_f(p^\ell)\lambda_f(p^{\ell+2})}{p^{\ell(1+\alpha+\beta+2s)}}  \right]\\
 & &\hspace{49mm} \times \prod_{p|N}\left[\left(1+\frac{\lambda_f(p)^2}{p^{1+u+v}}-\frac{\lambda_f(p)^2}{p^{1+\alpha+u+s}}-\frac{\lambda_f(p)^2}{p^{1+\beta+v+s}} \right)\sum_{\ell\geq 0} \frac{\lambda_f(p^\ell)^2}{p^{\ell(1+\alpha+\beta+2s)}} \right].
\end{eqnarray*}
\begin{normalsize}Thus, it comes from relations (\ref{eqn:relationmulti1}) and (\ref{eqn:relationmulti2}),\end{normalsize} 
\begin{eqnarray*}
& & \mathcal{P} =L(f\times f, 1+\alpha+\beta+2s) \prod_{p|N}\left[1+\frac{\lambda_f(p)^2}{p^{1+u+v}}-\frac{\lambda_f(p)^2}{p^{1+\alpha+u+s}}-\frac{\lambda_f(p)^2}{p^{1+\beta+v+s}} \right] \\
& &\hspace{66mm}\times  \prod_{p\nmid N} \left[   \left(1+\frac{\lambda_f(p)^2}{p^{1+u+v}}+\frac{1}{p^{2(1+u+v)}}\right)\left( 1-\frac{1}{p^{2(1+\alpha+\beta+2s)}} \right) \right. \\
& & \hspace{46mm}-\lambda_f(p)^2\left( 1-\frac{1}{p^{1+\alpha+\beta+2s}} \right)\left(\frac{1}{p^{1+v+\beta+s}} 
 + \frac{1}{p^{1+u+\alpha+s}}\right)\left(1+\frac{1}{p^{1+u+v}}\right)\\
 & & \hspace{33mm} + \left.\left( 1-\frac{1}{p^{1+\alpha+\beta+2s}} \right)\left( \lambda_f(p^2)-\frac{1}{p^{1+\alpha+\beta+2s}} \right)\left(\frac{1}{p^{2(1+v+\beta+s)}}+\frac{1}{p^{2(1+u+\alpha+s)}}\right)\right] \\
 & & \hspace{4mm}= L(f\times f, 1+\alpha+\beta+2s) \prod_{p}\left[1+\frac{\lambda_f(p)^2}{p^{1+u+v}}-\frac{\lambda_f(p)^2}{p^{1+\alpha+u+s}}-\frac{\lambda_f(p)^2}{p^{1+\beta+v+s}}+\chi_0(p)E_p \right] 
 \end{eqnarray*}
\begin{normalsize}where\end{normalsize} 
 \begin{align*}
 E_p = \hspace{155mm} \\ \frac{1}{p^2}\left[\frac{1}{p^{2(u+v)}} -\frac{1}{p^{2(\alpha +\beta +2s)}} -\frac{\lambda_f(p)^2}{p^s}\left(\frac{1}{p^{u+\alpha}}+\frac{1}{p^{v+\beta }} \right)\left(\frac{1}{p^{u+v}}-\frac{1}{p^{\alpha + \beta +2 s }} \right)+\frac{\lambda_f(p^2)}{p^{2s}}\left( \frac{1}{p^{2(u+\alpha ) }} + \frac{1}{p^{2(v+\beta ) }} \right) \right]\\
 +\frac{\lambda_f(p)^2}{p^{3+\alpha +\beta +3s }}\left[\frac{1}{p^{u+v }}\left(\frac{1}{p^{u+\alpha  }} + \frac{1}{p^{v+\beta  }} \right) - \frac{1}{p^{s}}\left( \frac{1}{p^{2(u+\alpha ) }} + \frac{1}{p^{2(v+\beta ) }}\right)-\frac{1}{p^{u+v+\alpha+\beta+s}} \right] \\
 +\frac{1}{p^{4+2(\alpha +\beta +2s) }}\left[\frac{1}{p^{2(u+\alpha+s)}}+\frac{1}{p^{2(v+\beta+s)}} -\frac{1}{p^{2(u+v)}} \right].
 \end{align*}
\begin{normalsize}Since the Rankin-Selberg $L$-function $L(f\times f,z)$ admits, for $\Re z > 1$, the Euler product\end{normalsize}
 \begin{eqnarray*} L(f\times f,z) =  \prod_{p\in \mathcal{P}} L_p(f\times f,z)  \end{eqnarray*}
\begin{normalsize}with \end{normalsize}
  \begin{eqnarray*}
  L_p(f\times f,z) & = & \left( 1-\frac{\alpha_f(p)^2 }{p^z} \right)^{-1}\left( 1-\frac{\alpha_f(p)\beta_f(p) }{p^z} \right)^{-2}\left( 1-\frac{\beta_f(p)^2 }{p^z} \right)^{-1} \\
   & = &\left(1-\frac{\lambda_f(p)^2}{p^z}+\chi_0(p)\left[\frac{2+\lambda_f(p)^2}{p^{2z}} -\frac{\lambda_f(p)^2}{p^{3z}}+\frac{1}{p^{4z}}\right] \right)^{-1} ,
 \end{eqnarray*}
\begin{normalsize}we may write\end{normalsize}
\begin{align*}
 1+\frac{\lambda_f(p)^2}{p^{1+u+v}}-\frac{\lambda_f(p)^2}{p^{1+\alpha+u+s}}-\frac{\lambda_f(p)^2}{p^{1+\beta+v+s}}+\chi_0(p)E_p  = \frac{ L_p(f\times f,1+u+v)}{ L_p(f\times f,1+u+\alpha +s) L_p(f\times f,1+v+\beta +s)} \\ \times \left[1+L_p(f\times f,1+u+\alpha +s) L_p(f\times f,1+v+\beta +s)\sum_{r= 2}^8\sum_{\ell}\frac{a_{r,\ell}(p)}{p^{r+X_{r,\ell}(u,v,\alpha,\beta,s)}} \right]
  \end{align*}
\begin{normalsize}where the sum over $\ell$ is finite, $X_{r,\ell}$ are linear forms in $u,v,\alpha,\beta,s$ and $a_{r,\ell}(p)$ are  complex numbers with $|a_{r,\ell}(p)|\ll 1$. As a result, we obtain\end{normalsize}
 \begin{eqnarray*}
 \mathcal{P}=\frac{L(f\times f, 1+\alpha+\beta+2s)L(f\times f, 1+u+v)}{L(f\times f, 1+\alpha+u+s)L(f\times f, 1+\beta+v+s)}A_{\alpha,\beta}(u,v,s) 
\end{eqnarray*}
\begin{normalsize}where \end{normalsize}
\begin{eqnarray*}
A_{\alpha,\beta}(u,v,s) = \prod_{p\in \mathcal{P}}\left[1+ \sum_{r,\ell} O\left( \frac{1}{p^{r+X_{r,\ell}(\Re u,\Re v,\Re \alpha,\Re \beta,\Re s)}}\right) \right].
\end{eqnarray*}
\end{small}Then  $A_{\alpha,\beta}(u,v,s)$ is  an absolutely convergent Euler product in $\{\Re( \alpha+u+s)>-1\}\cap \{\Re(\beta+v+s)>-1\} \cap_{r,\ell} \{X_{r,\ell}(\Re u,\Re v,\Re \alpha,\Re \beta,\Re s)>1-r \}$. Making explicit all the linear forms $X_{2,\ell}$, we obtain
$ \{X_{2,\ell}(\Re u,\Re v,\Re \alpha,\Re \beta,\Re s)>-1 \}=\Omega_{\alpha,\beta}$.
Similarly, writing all the linear forms $X_{r,\ell}$, we prove that for $(u,v,s)$ in $\Omega_{\alpha,\beta}$ then $X_{r,\ell}(\Re u,\Re v,\Re \alpha,\Re \beta,s) > -r/2\geq  1-r$. As a result, $A_{\alpha,\beta}(u,v,s)$ is an absolutely convergent Euler product on $\Omega_{\alpha,\beta}$ and defines an holomorphic function on $\Omega_{\alpha,\beta}$.
\begin{flushright}
$\Box$
\end{flushright} 

\begin{lemma} \label{lem:produiteulerienen0} We have 
$$ A_{0,0}(0,0,0)=1 . $$
\end{lemma}
Proof: Thanks to lemma \ref{eqn:expressionproduiteulerienenO}, if $\Re s>0$, we may write
\begin{eqnarray*} A_{0,0}(s,s,s)= \sum_{\substack{a,b,m,n \geq 1 \\am=bn }}\frac{\mu_f(a)\mu_f(b)\lambda_f(m)\lambda_f(n)}{(ambn)^{\frac12+s}} =\sum_{a,m\geq 1}\frac{\mu_f(a)\lambda_f(m)}{(am)^{1+2s}}\sum_{n |am}\mu_f(am/n)\lambda_f(n) . 
\end{eqnarray*}
Since $(\mu_f(n))$ is the convolution inverse of $(\lambda_f(n))$, we have 
$$\sum_{n |d}\mu_f(d/n)\lambda_f(n) =\delta(d).  $$
Thus
$$A_{0,0}(s,s,s) = \mu_f(1)\lambda_f(1) = 1 . $$
To conclude, we extend this relation to $s=0$ by continuity in the half plane $\Re s>0$.
\begin{flushright}
$\Box$
\end{flushright}
\begin{lemma} \label{lem:IPPEMcL} For all non negative integer $a$, we have 
$$ \sum_{n \leq M }\frac{\lambda_f(n)^2}{n}\left(\ln\frac{M}{n} \right)^a =\frac{\underset{s=1}{Res\,}L(f\times f,s)}{\zeta^{(N)}(2)}\int_{1}^{M}\frac{1}{r}\left( \ln  \frac{M}{r} \right)^a dr + O \left((\ln M)^a  \right). $$
\end{lemma}
Proof: We may find in \cite{Ra} the following asymptotic behaviour $$\sum_{n\leq x} \lambda_f(n)^2 = x\frac{\underset{s=1}{Res\,}L(f\times f,s)}{\zeta^{(N)}(2)} +O(x^{3/5}).$$
That is why, after one integration by parts, we get
$$ \sum_{n \leq M }\frac{\lambda_f(n)^2}{n}\left(\ln\frac{M}{n} \right)^a =\frac{\underset{s=1}{Res\,}L(f\times f,s)}{\zeta^{(N)}(2)}\sum_{n\leq M } \frac{1}{n}\left(\ln\frac{M}{n} \right)^a + O \left((\ln M)^a  \right) . $$
Finally, the Euler Maclaurin formula gives us $$ \sum_{n\leq M } \frac{1}{n}\left(\ln\frac{M}{n} \right)^a  = \int_{1}^{M}\frac{1}{r}\left( \ln  \frac{M}{r} \right)^a dr + O \left((\ln M)^a  \right).  $$
\begin{flushright}
$\Box$
\end{flushright}

\subsubsection{Estimation of $I_{f}^{D_1}(\alpha, \beta)$}
For positive integers $i$ and $j$, and for any positive real number $\delta$, let
\begin{eqnarray}
J_{\alpha,\beta}(i,j) = \frac{1}{(2i\pi)^2}\int_{(\delta)}\int_{(\delta)}M^{u+v}\frac{L(f\times f, 1+u+v)A_{\alpha,\beta}(u,v,0)}{L(f\times f, 1+\alpha+u)L(f\times f, 1+\beta+v)}\frac{du}{u^{i+1}}\frac{dv}{v^{j+1}}.
\end{eqnarray}
\begin{lemma} \label{lem:lemmesurJ} Let $\alpha, \beta \ll 1/L$ be complex numbers. For positive integers $i$ and $j$, we have
\begin{eqnarray*}J_{\alpha,\beta}(i,j) = \frac{(\ln M)^{i+j-1}}{i!j!\underset{s=1}{Res\,}L(f\times f,s)} \left.\frac{d^2}{dxdy}\left[M^{\alpha x+\beta y}\int_{0}^{1}(x+u)^i(y+u)^jdu  \right]\right|_{x=y=0}\hspace{23mm}\\+O\left(L^{i+j-2}\left( 1+\frac{(\ln L)^2}{L^{i-1}}\right)\left( 1+\frac{(\ln L)^2}{L^{j-1}}\right)\right) .
\end{eqnarray*}
\end{lemma}
Proof: We use the Dirichlet series of $L(f\times f,s)$. Moving far to the right either  $u$ or $v$ integration line, we obtain
\begin{small}
$$ J_{\alpha,\beta}(i,j) = \sum_{n\leq M}\frac{\lambda_f(n)^2}{n} \frac{1}{(2i\pi)^2}\int_{(\delta)}\int_{(\delta)}\left(\frac{M}{n}\right)^{u+v}\frac{\zeta^{(N)}(2( 1+u+v))A_{\alpha,\beta}(u,v,0)}{L(f\times f, 1+\alpha+u)L(f\times f, 1+\beta+v)}\frac{du}{u^{i+1}}\frac{dv}{v^{j+1}} .$$ 
\begin{normalsize}Let \end{normalsize}
$$r_{\alpha,\beta}^{(i,j)}(u,v)= \left(\frac{M}{n}\right)^{u+v}\frac{\zeta^{(N)}(2( 1+u+v))A_{\alpha,\beta}(u,v,0)}{L(f\times f, 1+\alpha+u)L(f\times f, 1+\beta+v)}\frac{1}{u^{i+1}v^{j+1}} . $$\end{small}

\noindent We consider the contour $\gamma$ defined by $\gamma = \gamma_1 \cup \gamma_2 \cup \gamma_3$ with $c>0$, $Y\gg 1$ large and
\begin{eqnarray*}
\gamma_1 = \{ i\tau , \;|\tau|\geq Y \},\hspace{2mm} \gamma_2 = \left\{ \sigma \pm i Y, \;\frac{-c}{\ln Y} \leq \sigma \leq 0 \right\} \mbox{ and } \gamma_3 = \left\{ \frac{-c}{\ln Y}+ i\tau, \; |\tau|\leq Y \right\}.
\end{eqnarray*}
By the standard zero-free region of $L(f\times f,s)$, we move the integration lines $\Re u = \Re v=\delta$ to $u$ and $v$ on $\gamma$. Thus, 
\begin{small}
\begin{align*}
 & \frac{1}{(2i\pi)^2}\int_{(\delta)}\int_{(\delta)} r_{\alpha,\beta}^{(i,j)}(u,v)dudv  \\  = & \underset{u=0}{Res}\frac{1}{2i\pi} \int_{\Re v = \delta}r_{\alpha,\beta}^{(i,j)}(u,v)dv +\frac{1}{(2i\pi)^2}\int_{u \in \gamma }\int_{\Re v =\delta} r_{\alpha,\beta}^{(i,j)}(u,v)dvdu \\
 = &\underset{u=v=0}{Res}\, r_{\alpha,\beta}^{(i,j)}(u,v) +\underset{u=0}{Res}\frac{1}{2i\pi} \int_{v \in \gamma}r_{\alpha,\beta}^{(i,j)}(u,v)dv + \underset{v=0}{Res}\frac{1}{2i\pi} \int_{u \in \gamma}r_{\alpha,\beta}^{(i,j)}(u,v)du + \frac{1}{(2i\pi)^2}\int_{\gamma}\int_{\gamma} r_{\alpha,\beta}^{(i,j)}(u,v)dudv .
\end{align*}
$\bullet$ \begin{normalsize}We begin with the estimation of\end{normalsize} $\underset{u=0}{Res}\frac{1}{2i\pi} \int_{v \in \gamma}r_{\alpha,\beta}^{(i,j)}(u,v)dv$. \begin{normalsize}We express the residue as a contour integral on a circle with radius $1/L$. We get
\end{normalsize}
\begin{eqnarray*}
& & \underset{u=0}{Res}\frac{1}{2i\pi} \int_{v \in \gamma}r_{\alpha,\beta}^{(i,j)}(u,v)dv  \\ & = & \frac{1}{(2i\pi)^2}\int_{v\in \gamma} \frac{\left(\frac{M}{n}\right)^{v}}{L(f\times f, 1+\beta+v)}\oint_{D(0,L^{-1})} \left(\frac{M}{n}\right)^{u}\frac{\zeta^{(N)}(2( 1+u+v))A_{\alpha,\beta}(u,v,0)}{L(f\times f, 1+\alpha+u)}\frac{du}{u^{i+1}} \frac{dv}{v^{j+1}}.
\end{eqnarray*} 
\begin{normalsize}Furthermore, since $| \zeta^{(N)}(2( 1+u+v))A_{\alpha,\beta}(u,v,0) |\ll 1$ in our range of integration and since
\end{normalsize}
$$\frac{1}{L(f\times f, 1+\alpha+u)} \ll \alpha+u \ll L^{-1} \mbox{ because } u \asymp 1/L,$$
\begin{normalsize}we obtain 
\end{normalsize}
\begin{eqnarray*}
\underset{u=0}{Res}\frac{1}{2i\pi} \int_{v \in \gamma}r_{\alpha,\beta}^{(i,j)}(u,v)dv \ll L^{i-1} \int_{v\in \gamma} \frac{\left(\frac{M}{n}\right)^{\Re v}}{\left|L(f\times f, 1+\beta+v)\right|}\frac{dv}{|v|^{j+1}}.
\end{eqnarray*}
\begin{normalsize}As we have $-\frac{L'(f\times f,z)}{L(f\times f, z)} = \sum_{n\geq 1}\frac{\Lambda_f(n)}{n^z}$ with $\Lambda_f(n)\geq 0$, we deduce (confer  \cite{Te} part 3.10 and \cite{IK} part 5.3) that
\end{normalsize}
$$\frac{1}{L(f\times f, \sigma+i \tau )} \ll \ln |\tau|, $$
\begin{normalsize}it comes
\end{normalsize}
\begin{eqnarray}
& & \int_{v \in \gamma} \frac{\left(\frac{M}{n}\right)^{\Re v}}{\left|L(f\times f, 1+\beta+v)\right|}\frac{dv}{|v|^{j+1}} \nonumber \\ & \ll & \int_{|\tau|\geq Y}\frac{\ln \tau}{|\tau|^{j+1}}d\tau + \ln Y \int^0_{-\frac{c}{\ln Y}}\frac{d\sigma}{|\sigma+iY|^{j+1}}  
+\left(\frac{M}{n}\right)^{-\frac{c}{\ln Y}}\ln Y \int_{|\tau|\leq Y} \frac{d\tau}{\left| \tau  - i\frac{c}{\ln Y}\right|^{j+1}} \nonumber \\
 &  \ll & \frac{\ln Y}{Y^j}+(\ln Y)^{j+1}\left(\frac{M}{n}\right)^{-\frac{c}{\ln Y}}. \label{eqn:estimationdansregionsanszero}
\end{eqnarray}
\begin{normalsize}As a result, we get the bound\end{normalsize}
\begin{eqnarray}\underset{u=0}{Res}\frac{1}{2i\pi} \int_{v \in \gamma}r_{\alpha,\beta}^{(i,j)}(u,v)dv \ll L^{i-1}\ln Y\left( \frac{1}{Y^j}+(\ln Y)^{j}\left(\frac{M}{n}\right)^{-\frac{c}{\ln Y}} \right).\label{eqn:estimationdansregionsanszero1}\end{eqnarray} 
$\bullet$ \begin{normalsize}Since\end{normalsize} $r_{\alpha,\beta}^{(i,j)}(u,v) = r_{\beta,\alpha}^{(j,i)}(v,u)$, \begin{normalsize}it comes immediately from the previous bound\end{normalsize} 
\begin{eqnarray}\underset{v=0}{Res}\frac{1}{2i\pi} \int_{u \in \gamma}r_{\alpha,\beta}^{(i,j)}(u,v)du \ll L^{j-1}\ln Y\left( \frac{1}{Y^i}+(\ln Y)^{i}\left(\frac{M}{n}\right)^{-\frac{c}{\ln Y}} \right).\label{eqn:estimationdansregionsanszero2}\end{eqnarray} 
$\bullet$ \begin{normalsize}By the relation (\ref{eqn:estimationdansregionsanszero})\end{normalsize}, we bound $\int_{\gamma}\int_{\gamma} r_{\alpha,\beta}^{(i,j)}(u,v)dudv$.
\begin{align}\frac{1}{(2i\pi)^2}\int_{\gamma}\int_{\gamma} r_{\alpha,\beta}^{(i,j)}(u,v)dudv&=\frac{1}{(2i\pi)^2}\int_{\gamma}\int_{\gamma} \left(\frac{M}{n}\right)^{u+v}\frac{\zeta^{(N)}(2( 1+u+v))A_{\alpha,\beta}(u,v,0)}{L(f\times f, 1+\alpha+u)L(f\times f, 1+\beta+v)}\frac{1}{u^{i+1}v^{j+1}} \nonumber \\
&\ll \int_\gamma \frac{\left(\frac{M}{n}\right)^{\Re v}}{\left|L(f\times f, 1+\beta+v)\right|}\frac{dv}{|v|^{j+1}}\int_\gamma \frac{\left(\frac{M}{n}\right)^{\Re u}}{\left|L(f\times f, 1+\alpha+u)\right|}\frac{du}{|u|^{i+1}} \nonumber \\
&\ll \left(\frac{\ln Y}{Y^j}+(\ln Y)^{j+1}\left(\frac{M}{n}\right)^{-\frac{c}{\ln Y}}\right)\left(\frac{\ln Y}{Y^i}+(\ln Y)^{i+1}\left(\frac{M}{n}\right)^{-\frac{c}{\ln Y}} \right)\nonumber \\ &\ll (\ln Y)^2 \left( \frac{1}{Y^{i+j}}+(\ln Y)^{i+j}\left(\frac{M}{n}\right)^{-\frac{c}{\ln Y}}\right). \label{eqn:estimationdansregionsanszero3}
\end{align}
$\bullet$ \begin{normalsize}In addition, for any positive integer $\ell$, let\end{normalsize}
$$W_\ell=\sum_{n\leq M}\frac{\lambda_f(n)^2}{n} \left( \frac{1}{Y^\ell}+(\ln Y)^{\ell}\left(\frac{M}{n}\right)^{-\frac{c}{\ln Y}}\right).$$
\begin{normalsize} We can bound\end{normalsize}
\begin{align*}
W_\ell & \ll \frac{1}{Y^\ell}\sum_{n\leq \frac{M}{(Y\ln Y)^{\ell(\ln Y)/c}}}\frac{\lambda_f(n)^2}{n} + (\ln Y)^{\ell}\sum_{ \frac{M}{(Y\ln Y)^{\ell(\ln Y)/c}} \leq n\leq M}\frac{\lambda_f(n)^2}{n}\left(\frac{M}{n}\right)^{-\frac{c}{\ln Y}} \\
& \ll  \frac{L}{Y^\ell}+ (\ln Y)^{\ell} \sum_{0\leq d \leq  \frac{\ln(Y\ln Y)}{\ln 2}} 2^{-d \ell}\sum_{\frac{M}{\left(2^{d +1}\right)^{\ell(\ln Y)/c}}\leq n \leq \frac{M}{\left(2^{d}\right)^{\ell(\ln Y)/c}}}\frac{\lambda_f(n)^2}{n} \\ & \ll  \frac{L}{Y^\ell} + (\ln Y)^{\ell+1}. 
\end{align*}\begin{normalsize}
As a consequence, by relations (\ref{eqn:estimationdansregionsanszero1}), (\ref{eqn:estimationdansregionsanszero2}) and (\ref{eqn:estimationdansregionsanszero3}), we obtain\end{normalsize} 
\begin{align*}J_{\alpha,\beta}(i,j)& =   \sum_{n\leq M}\left[ \frac{\lambda_f(n)^2}{n}\underset{u=v=0}{Res}\, r_{\alpha,\beta}^{(i,j)}(u,v)\right] +O\left(  L^{i-1}W_j\ln Y  + L^{j-1}W_i\ln Y + W_{i+j}(\ln Y)^2 \right) \\
 & =  \sum_{n\leq M} \left[ \frac{\lambda_f(n)^2}{n}\underset{u=v=0}{Res}\, r_{\alpha,\beta}^{(i,j)}(u,v)\right] +O\left[ L^{i-1}\ln Y\left(\frac{L}{Y^j}+(\ln Y)^{j} \right)+L^{j-1}\ln Y\left(\frac{L}{Y^i}+(\ln Y)^{i} \right)\right. \\
&  \hspace{96mm}\left. +(\ln Y)^2 \left( \frac{L}{Y^{i+j}}+(\ln Y)^{i+j} \right)\right]. \end{align*}
\begin{normalsize}Choosing $Y=L$, which is allowed, the error term becomes $O\left(\left(L^{i-1}+(\ln L)^{2}\right)\left(L^{j-1}+(\ln L)^{2}\right)\right)$. Then\end{normalsize}
\begin{eqnarray}
J_{\alpha,\beta}(i,j) = \sum_{n\leq M}\left[\frac{\lambda_f(n)^2}{n}\underset{u=v=0}{Res}\, r_{\alpha,\beta}^{(i,j)}(u,v)\right]+O\left(L^{i+j-2}\left(1+\frac{(\ln L)^2}{L^{i-1}} \right)\left(1+\frac{(\ln L)^2}{L^{j-1}} \right) \right) \label{eqn:estiJetape1}
\end{eqnarray}
$\bullet$ \begin{normalsize}We finish with the estimation of $\underset{u=v=0}{Res}\, r_{\alpha,\beta}^{(i,j)}(u,v)$. To do so again, we express the residue as a contour integral on a circle with radius $1/L$. Thus, 
\end{normalsize}
\begin{align*} \underset{u=v=0}{Res}\, r_{\alpha,\beta}^{(i,j)}(u,v)= \ \frac{1}{(2i\pi)^2}\oint_{D(0,L^{-1})}\oint_{D(0,L^{-1})}\left(\frac{M}{n}\right)^{u+v}\frac{\zeta^{(N)}(2( 1+u+v))A_{\alpha,\beta}(u,v,0)}{L(f\times f, 1+\alpha+u)L(f\times f, 1+\beta+v)}\frac{du}{u^{i+1}}\frac{dv}{v^{j+1}} 
\end{align*}
\begin{normalsize}Furthermore, since $ u \asymp  v \asymp 1/L$, we have 
\end{normalsize}
\begin{eqnarray*}
\zeta^{(N)}(2( 1+u+v)) & = &\zeta^{(N)}(2)+O(1/L),\\  A_{\alpha,\beta}(u,v,0) & = & A_{0,0}(0,0,0)+O(1/L), \\
 \frac{1}{L(f\times f, 1+\alpha+u)} & = & \frac{\alpha+u}{\underset{s=1}{Res\,}L(f\times f,s)}(1+O(1/L)),\\  \frac{1}{L(f\times f, 1+\beta+v)} & = & \frac{\beta+v}{\underset{s=1}{Res\,}L(f\times f,s)}(1+O(1/L)).
\end{eqnarray*}
\begin{normalsize}With lemma \ref{lem:produiteulerienen0}, we obtain
\end{normalsize}
$$\frac{\zeta^{(N)}(2( 1+u+v))A_{\alpha,\beta}(u,v,0)}{L(f\times f, 1+\alpha+u)L(f\times f, 1+\beta+v)}=(\alpha+u)(\beta+v)\frac{\zeta^{(N)}(2)}{[\underset{s=1}{Res\,}L(f\times f,s)]^2}+O(1/L^3). $$
\begin{normalsize}Then 
\end{normalsize}
\begin{align*} \underset{u=v=0}{Res}\, r_{\alpha,\beta}^{(i,j)}(u,v) = \frac{\zeta^{(N)}(2)}{[\underset{s=1}{Res\,}L(f\times f,s)]^2}\left[\frac{1}{2i\pi}\oint_{D(0,L^{-1})} \left(\frac{M}{n}\right)^{u}\frac{\alpha+u}{u^{i+1}}du \right]\left[ \frac{1}{2i\pi}\oint_{D(0,L^{-1})} \left(\frac{M}{n}\right)^{v}\frac{\beta+v}{v^{j+1}}dv\right]\\
+O\left(L^{i+j-3}\right).
\end{align*}
\begin{normalsize}In addition, thanks to  Cauchy formula, we compute for any positive integer $\ell$
\end{normalsize}
\begin{eqnarray*}\frac{1}{2i\pi}\oint_{D(0,L^{-1})} \left(\frac{M}{n}\right)^{u}\frac{\alpha+u}{u^{\ell+1}}du  & = & \left. \frac{d}{dx}\left[ \frac{e^{\alpha x}}{2i\pi}\oint_{D(0,L^{-1})}\left(\frac{M}{n}e^x\right)^{u}\frac{du}{u^{\ell+1}}\right]\right|_{x=0}
\\& =& \left. \frac{d}{dx}\left[ e^{\alpha x} \frac{1}{\ell !}\left. \frac{d^\ell}{du^\ell} \left[  \left(\frac{M}{n}e^x\right)^{u}\right]\right|_{u=0} \right]\right|_{x=0} \\ & =& \frac{1}{\ell !}\left. \frac{d}{dx}\left[ e^{\alpha x}\left(x+\ln \frac{M}{n}\right)^\ell \right]\right|_{x=0}\\& = &\frac{1}{\ell !} \left[\alpha  \left(\ln \frac{M}{n}\right)^\ell +\ell \left(\ln \frac{M}{n}\right)^{\ell -1} \right].
\end{eqnarray*}
\begin{normalsize}Then, we can write
\end{normalsize}\begin{align*}
\underset{u=v=0}{Res}\, r_{\alpha,\beta}^{(i,j)}(u,v) = \frac{\zeta^{(N)}(2)}{i!j![\underset{s=1}{Res\,}L(f\times f,s)]^2}\left[\alpha\beta \left(\ln\frac{M}{n}\right)^{i+j} +(\alpha i+\beta j)\left(\ln\frac{M}{n}\right)^{i+j-1}  + ij \left(\ln\frac{M}{n}\right)^{i+j-2} \right] \\ + O\left(L^{i+j-3}\right).
\end{align*}\begin{normalsize}From (\ref{eqn:estiJetape1}) and by lemma \ref{lem:IPPEMcL},  we get
\end{normalsize}
\begin{align*}
J_{\alpha, \beta}(i, j ) = \frac{[\underset{s=1}{Res\,}L(f\times f,s)]^{-1}}{i!j!}\int_1^M \left( \alpha\left(\ln\frac{M}{r} \right)^{i}+i\left(\ln\frac{M}{r} \right)^{i-1}\right)\left( \beta\left(\ln\frac{M}{r} \right)^{j}+j\left(\ln\frac{M}{r} \right)^{j-1}\right)\frac{dr}{r} \\ + O\left(L^{i+j-2}\left(1+\frac{(\ln L)^2}{L^{i-1}} \right)\left(1+\frac{(\ln L)^2}{L^{j-1}} \right) \right).
\end{align*}
\begin{normalsize}Changing the variable $r$ to $u$ with $r=M^{1-u}$ concludes the proof.
\end{normalsize}
\end{small}
\begin{flushright}
$\Box$
\end{flushright}

\begin{lemma} \label{pro:termediagpartie1} If  $0<\nu<1$ and $\alpha, \beta \ll L^{-1}$ are complex numbers with $|\alpha+\beta | \gg L^{-1}$, then
$$I_{f}^{D_1}(\alpha, \beta ) = \frac{\widehat{w}(0)}{(\alpha+\beta)\ln M} \left.\frac{d^2}{dxdy}\left[M^{\alpha x+\beta y}\int_{0}^{1}P(x+u)P(y+u)du  \right]\right|_{x=y=0}+O\left( \frac{T(\ln L)^4}{L}\right) .$$
\end{lemma}
Proof: We use the Mellin transformation to write
$$ \left(\frac{\ln M/a}{\ln M}  \right)^i = \left\{ \begin{array}{l} \frac{i!}{(\ln M)^i}\frac{1}{2i\pi}\int_{(1)}\left(\frac{M}{a}\right)^v\frac{dv}{v^{i+1}}  \hspace{2mm}\mbox { if }1\leq a \leq M , \\ 
0 \hspace{3mm} \mbox{ if } a>M . \end{array} \right.  $$
Set $P(X)=\sum_{i=1}^{\deg P}a_iX^i$. Thus, from  relation (\ref{eqn:defIfD1}) which defines $I_f^{D_1}(\alpha, \beta )$, we get
\begin{multline*}I_f^{D_1}(\alpha, \beta ) = \int_{-\infty}^{+\infty}w(t)\sum_{i,j\geq 1}\frac{a_ia_ji!j!}{(\ln M)^{i+j}} \frac{1}{(2i\pi)^3}\int_{(1)}\int_{(1)}\int_{(1)}\frac{M^{u+v}G(s)}{s}g_{\alpha,\beta}(s,t) \\
  \sum_{\substack{a,b,m,n\geq 1 \\am=bn}}\frac{\mu_f(a)\mu_f(b)\lambda_f(m)\lambda_f(n)}{a^{\frac12+v}b^{\frac12+u}m^{\frac12+\alpha+s}n^{\frac12+\beta+s}} ds \frac{du}{u^{j+1}}\frac{dv}{v^{i+1}}dt  .
 \end{multline*}
Due to lemma \ref{eqn:expressionproduiteulerienenO}, we can write
\begin{multline*}
I_f^{D_1}(\alpha, \beta) = \int_{-\infty}^{+\infty}w(t)\sum_{i,j\geq 1}\frac{a_ia_ji!j!}{(\ln M)^{i+j}} \frac{1}{(2i\pi)^3}\int_{(1)}\int_{(1)}\int_{(1)}\frac{M^{u+v}G(s)}{s}g_{\alpha,\beta}(s,t)A_{\alpha,\beta}(u,v,s) \\
 \frac{L(f\times f, 1+\alpha+\beta+2s)L(f\times f, 1+u+v)}{L(f\times f, 1+\alpha+u+s)L(f\times f, 1+\beta+v+s)} ds \frac{du}{u^{j+1}}\frac{dv}{v^{i+1}}dt .
\end{multline*}
We specialize $G$ in $$G(s)= e^{s^2}\frac{(\alpha+\beta)^2-(2s)^2}{(\alpha+\beta)^2}.$$
As a result, $G(s)L(f\times f, 1+\alpha+\beta+2s)$ is an entire function. First, we move the integration lines from $\Re u=$Re $\Re v=1$ to $\Re u=\Re v=\delta$, with $\delta$ small, in order to assure the absolute convergence of $A_{\alpha,\beta}(u,v,s)$. Secondly, we move the integration line  from $\Re s=1$ to $\Re s=-\delta+\epsilon$, with $0<\epsilon<\delta$, crossing a pole at $s=0$. Since $t\asymp T$, $\nu <1$ and $g_{\alpha, \beta}(s,t) \ll T^{2s}$, we can bound
\begin{multline*}\int_{-\infty}^{+\infty}w(t)\sum_{i,j}\frac{a_ia_ji!j!}{(\ln M)^{i+j}} \frac{1}{(2i\pi)^3}\int_{Re \,u=\delta}\int_{Re \,v=\delta}\int_{Re \, s =-\delta+\epsilon}\frac{M^{u+v}G(s)}{s}g_{\alpha,\beta}(s,t)A_{\alpha,\beta}(u,v,s) \\
 \frac{L(f\times f, 1+\alpha+\beta+2s)L(f\times f, 1+u+v)}{L(f\times f, 1+\alpha+u+s)L(f\times f, 1+\beta+v+s)} ds \frac{du}{u^{j+1}}\frac{dv}{v^{i+1}}dt  \\
\ll \int_{-\infty}^{+\infty}|w(t)|dt\, T^{2(-\delta+\epsilon)}M^{2\delta}  \ll T^{1-(2-2\nu)\delta+\epsilon}\ll T^{1-\varepsilon}.
\end{multline*}
for sufficiently small $\varepsilon $. Then, using some previous notations, this estimation gives
\begin{small}$$ I_f^{D_1}(\alpha, \beta) = \widehat{w}(0)L(f\times f, 1+\alpha +\beta)\sum_{i,j}\frac{a_ia_ji!j!}{(\ln M)^{i+j}}J_{\alpha, \beta}(i,j)+O\left(T^{1-\varepsilon}\right).  $$\end{small}Thanks to lemma \ref{lem:lemmesurJ} and since \begin{small}$$ L(f \times f, 1+\alpha+\beta) = \frac{\underset{s=1}{Res\,}L(f\times f,s)}{\alpha+\beta} +O(1), $$\end{small}we get
\begin{small}
\begin{multline*}
I_f^{D_1}(\alpha, \beta) = \widehat{w}(0)  \left[\frac{\underset{s=1}{Res\,}L(f\times f,s)}{\alpha+\beta} +O(1)  \right] \\ \times  \left[ \frac{1}{\underset{s=1}{Res\,}L(f\times f,s)\ln M} \left.\frac{d^2}{dxdy}\left[M^{\alpha x+\beta y}\int_{0}^{1} P(x+u)P(y+u)du \right] \right|_{x=y=0}+O\left(\frac{(\ln L)^4}{L^2}\right) \right] +O\left(T^{1-\varepsilon}\right) \\
 =    \frac{\widehat{w}(0)}{(\alpha+\beta)\ln M}\left.\frac{d^2}{dxdy}\left[M^{\alpha x+\beta y}\int_{0}^{1} P(x+u)P(y+u)du \right] \right|_{x=y=0}+O\left( \frac{T(\ln L)^4}{(\alpha+\beta)L^2} \right)+O(T/L).
\end{multline*} \end{small}We conclude using the assumption $|\alpha+\beta|\gg L^{-1}$.
\begin{flushright}
$\Box$
\end{flushright}
\begin{remark} Sometimes, it may be useful to consider the relation 
\begin{small}\begin{align}
\left.\frac{d^2}{dxdy}\left[M^{\alpha x+\beta y}\int_{0}^{1}P(x+u)P(y+u)du  \right]\right|_{x=y=0}  
 = \int_{0}^{1}(P'(u)+\alpha\ln M P(u))(P'(u)+\beta\ln M P(u))du  .\label{eqn:expressionbis}
\end{align}\end{small}
\end{remark}

\begin{lemma}\label{pro:termediagpartie2} If $0<\nu<1$ and $\alpha, \beta \ll L^{-1}$ are complex numbers with $|\alpha+\beta | \gg L^{-1}$, then
$$ I_{f}^{D2}(\alpha, \beta) =  T^{-2(\alpha+\beta)}I_{f}^{D1}(-\beta, -\alpha)+O(T/L) . $$
\end{lemma}
Proof: We write
\begin{small}
$$ I_{a,b}^{D_2}(\alpha,\beta) =  \sum_{am=bn}\frac{\lambda_f(m)\lambda_f(n)}{m^{\frac12-\beta}n^{\frac12-\alpha}}\int_{-\infty}^{+\infty}w(t)\left( \frac{t\sqrt{N}}{2\pi}\right)^{-2(\alpha+\beta)}V_{-\beta,-\alpha}(mn,t)dt+ O\left(\sum_{\substack{am=bn \\ mn \ll T^{2+\varepsilon} }} \frac{\left|\lambda_f(m)\lambda_f(n)\right|}{m^{\frac12-\beta}n^{\frac12-\alpha}}\right).$$
\end{small}Let $a'=a/(a,b)$ and $b'=b/(a,b)$. Then, for all $\delta>0$, the error term in the previous relation becomes
\begin{small}$$\sum_{\substack{ am=bn \\ mn \ll T^{2+\varepsilon}  } } \frac{\left|\lambda_f(m)\lambda_f(n)\right|}{m^{\frac12-\beta}n^{\frac12-\alpha}} \ll \frac{1}{(a'b')^{1/2-\delta}}\sum_{k \ll \frac{T^{1+\frac{\varepsilon}{2}}}{\sqrt{a'b'}}}\frac{1}{k^{1-\alpha-\beta-2\delta}} \ll \frac{T^{\varepsilon}}{\sqrt{a'b'}}$$\end{small}Therefore, if $w_2(t)= w(t)\left( \frac{t\sqrt{N}}{2\pi}\right)^{-2(\alpha+\beta)}$, we may write
\begin{small}
\begin{align*}I_{f}^{D_2}(\alpha, \beta) = \sum_{a,b\leq M } \frac{\mu_f(a)\mu_f(b)}{\sqrt{ab}}P\left(\frac{\ln M/a}{\ln M}  \right)P\left(\frac{\ln M/b}{\ln M}  \right) \left[\sum_{am=bn}\frac{\lambda_f(m)\lambda_f(n)}{m^{\frac12-\beta}n^{\frac12-\alpha}}\int_{-\infty}^{+\infty}w_2(t)V_{-\beta,-\alpha}(mn,t)dt \right] \\  +O\left(T^{\varepsilon}\sum_{a,b\leq M }  \frac{1}{\sqrt{ab}\sqrt{a'b'}} \right)  .  \end{align*}
\end{small}We also have\begin{small} $$\sum_{a,b\leq M }  \frac{1}{\sqrt{ab}\sqrt{a'b'}} \ll \sum_{k\leq M} \frac{1}{k} \sum_{a,b\leq \frac{M}{k}}\frac{1}{ab} \ll \sum_{k\leq M} \frac{1}{k}\left(\ln\frac{M}{k}\right)^2 \ll \left(\ln M \right)^3$$\end{small}and since $w_2$ satisfies (\ref{eqn:hypwa}), (\ref{eqn:hypwb}) and (\ref{eqn:hypwc}), up to changing $w$ by $w_2$ and $(\alpha,\beta)$ by $(-\beta,-\alpha)$, lemma \ref{pro:termediagpartie1} gives  
\begin{small}$$I_{f}^{D_2}(\alpha, \beta) = \frac{\widehat{w_2}(0)}{(-\alpha-\beta)\ln M} \left.\frac{d^2}{dxdy}\left[M^{-\alpha x-\beta y}\int_{0}^{1}P(x+u)P(y+u)du  \right]\right|_{x=y=0}+O\left(T(\ln L)^4/L\right).$$\end{small}To conclude, due to the support of $w$, we can write $\left( \frac{t\sqrt{N}}{2\pi}\right)^{-2(\alpha+\beta)} =T^{-2(\alpha+\beta)}+O(1/L)$ which gives $\widehat{w_2}(0)=T^{-2(\alpha+\beta)}\widehat{w}(0)+O(T/L)$.
\begin{flushright}
$\Box$
\end{flushright}

\subsubsection{Proof of proposition \ref{pro:termmediag}}

From relation (\ref{eqn:splitTD}) and using lemma \ref{pro:termediagpartie2}, we can write 
\begin{eqnarray*}
I_{f}^{D}(\alpha,\beta)  & = &  I_{f}^{D1}(\alpha, \beta)+T^{-2(\alpha+\beta)}I_{f}^{D1}(-\beta, -\alpha)+O(T(\ln L)^4/L) \\
 & = & I_{f}^{D1}(\alpha, \beta)+I_{f}^{D1}(-\beta, -\alpha) + I_{f}^{D1}(-\beta, -\alpha)\left[T^{-2(\alpha+\beta)}-1\right]+O(T(\ln L)^4/L).
\end{eqnarray*}
Finally, using relation (\ref{eqn:expressionbis}) and lemma \ref{pro:termediagpartie1}, we have 
$$I_{f}^{D1}(\alpha, \beta)+I_{f}^{D1}(-\beta, -\alpha)=\widehat{w}(0)+O(T(\ln L)^4/L).$$
Combining these relations, we get the result.
\begin{flushright}
$\Box$
\end{flushright}

\section{Effective proportion of zeros on the critical line} \label{sec:EFPZCL}
In this section, we prove corollary \ref{cor:proportion}. From theorem \ref{pro:momentdordre2ramollilisse}, we may deduce the following theorem about mollified second moment of $L(f,s)$ and its derivative.  
\begin{theorem} \label{th:mm2} Let $Q$ be a polynomial with complex coefficients satisfying $Q(0)=1$. Let 
\begin{small}
$$V(s)= Q\left(-\frac{1}{2\ln T} \frac{d}{ds} \right)(L(f,s)).  $$
\end{small}
Then, if $\nu < \frac{1-2\theta}{4+2\theta}$, we have
\begin{small}
$$ \frac{1}{T}\int_{1}^{T}|V\psi(\sigma_0+it)|^2dt = c(P,Q,2R,\nu/2)+o(1) $$
\end{small}
where 
\begin{small}
\begin{eqnarray*} c(P,Q,r,\xi)  & = & 1+ \frac{1}{\xi}\int_{0}^{1}\int_{0}^{1}e^{2rs}\left[ \left. \frac{d}{dx}\left(e^{r\xi x}Q(s+\xi x)P(x+u) \right)\right|_{x=0} \right]^2du ds  .
  \end{eqnarray*}
\end{small}
\end{theorem} 
We do not give the proof of this theorem, which is essentially the same than the one of theorem 1 in \cite{Yo}. We may refer to part 2 and 3 of \cite{Yo} to find more details.
Now, $Q$ refers to a polynomial with complex coefficients of the shape 
\begin{eqnarray} \label{eqn:formepolynome}
Q(x)=1+\sum_{n=1}^M i^{n+1}\lambda_n\left[(1-2x)^n-1 \right]
\end{eqnarray} 
where $M$ is a positive integer and $(\lambda_1,...,\lambda_M)$ belongs to $\R^M$. 
\paragraph*{}Let $N_f(T)$ (resp. $N_{f,0}(T)$) be the number of non-trivial zeros $\rho$ (resp. on the critical line) of $L(f,s)$ with $0<\Im (\rho)\leq T$ for $f$ a holomorphic primitive cusp form of even weight, square-free level and trivial character.
 \begin{proposition} For $f$ a holomorphic primitive cusp form of even weight, square-free level and trivial character, $Q$ as in \emph{(\ref{eqn:formepolynome})}, we have
$$   \liminf_{T\rightarrow +\infty}\frac{N_{f,0}(T)}{N_f(T)} \geq \limsup_{T\rightarrow +\infty}\left[  1-\frac{1}{2R}\ln\left( \frac{1}{T}\int_{1}^{T}|V\psi(\sigma_0+it)|^2 dt \right)\right]$$
\end{proposition}
Therefore, using theorem \ref{th:mm2}, we get
$$ \liminf_{T\rightarrow +\infty}\frac{N_{f,0}(T)}{N_f(T)} \geq 1- \inf_{P,Q,R} \frac{1}{R}\ln c(P,Q,R,\nu/2). $$
The work of Kim and Sarnak (\cite{Ki}) gives $\theta=7/64$ then theorem \ref{th:mm2} gives $\nu=5/27$. The Ramanujuan-Petersson conjecture ($\theta=0$) gives $\nu=1/4$.
\begin{lemma}[\cite{C}, part 4] We have
$$\inf_{P} \frac{1}{R}\ln c(P,Q,R,\nu/2) = \frac{1}{R}\ln \left(\frac{1+|w(1)|^2}{2}+\frac{A\alpha}{\tanh \frac{\nu \alpha}{2} } \right)$$
where $w(x)=e^{Rx}Q(x)$, $A=\int_0^1 |w(x)|^2dx$, $B=\int_0^1 w(x)\overline{w'(x)}dx$, $C=\int_0^1 |w'(x)|^2dx$ and $\alpha= \frac{\sqrt{(B-\overline{B})^2+4AC}}{2A}$. 
\end{lemma}
For empirical reasons, we restrict ourselves to polynomials $Q$ with  real coefficients of the shape
\begin{eqnarray} Q(x)=1+\sum_{n=1}^N h_n\left[(1-2x)^{2n-1}-1 \right] \label{eqn:formepolynome2} \end{eqnarray}
where $N$ is a positive integer and $(h_1,...,h_N)$ belongs to $\R^N$.
Then  we get
\begin{eqnarray} \liminf_{T\rightarrow +\infty}\frac{N_{f,0}(T)}{N_f(T)} \geq 1 -\inf_{Q \mbox{\begin{scriptsize} real, \end{scriptsize}}R} \frac{1}{R}\ln \left(\frac{1+w(1)^2}{2}+\frac{\sqrt{AC}}{\tanh \left( \frac{\nu }{2}\sqrt{\frac{C}{A}}\right) } \right) . 
\end{eqnarray}
To obtain corollary \ref{cor:proportion}, we choose $N=4$, $R$ and $Q$ as in (\ref{eqn:formepolynome2}) where
 $R, h_1, h_2, h_3, h_4$ are given in the following table.
\begin{center}
\begin{tabular}{|c|c|c|c|}
\hline
 & $\nu=\frac{1}{6} \rule[-4mm]{0cm}{10mm} $ & $\nu=\frac{5}{27}$ & $\nu=\frac{1}{4}$ \\
\hline
$R$ & $6,6838894702116801322 $ &  $6,4278834168344993342$ & $5,6503610091685135131$\\
\hline
$h_1$ &  $1,6017785744634898860 $  & $1,5898336242677838745 $ & $ 1,5369390514358411982$\\
\hline
$h_2$&  $-3,0362512753510924917 $ & $ -2.8999828229132398066$  & $ -2,7929104872905007806$\\
\hline
$h_3$ &   $3,0757757634512927939 $ & $ 3.0171733454035522056$ & $ 2,7758193765120241770$\\
\hline
$h_4$ &  $-1,1407980564855935531 $ & $-1,1164150244992046552$& $ -1,0187870607687957034$\\
\hline 
\end{tabular}
\end{center}

\section{Non-mollified second integral moment}

This section contains the proof of theorem \ref{th:moment2nonramolli}. For more convenience, we set
$$M_{f,2}(\alpha,\beta)= \int_{-\infty}^{+\infty} w(t)L\left(f,\frac{1}{2}+\alpha+it\right)L\left(f,\frac{1}{2}+\beta-it\right)dt$$
where $w$ satisfies (\ref{eqn:hypwa}), (\ref{eqn:hypwb}) and (\ref{eqn:hypwc}). Applying the approximate equation (lemma \ref{lem:EFA}) and using some previous notations, we can write
$$M_{f,2}(\alpha,\beta)= I^{D_1}_{1,1}(\alpha,\beta)+I^{D_2}_{1,1}(\alpha,\beta)+I^{ND_1}_{1,1}(\alpha,\beta)+I^{ND_2}_{1,1}(\alpha,\beta). $$
Proposition \ref{pro:estiOD1} and corollary \ref{pro:estiOD2} allows to bound the off-diagonal contribution, then
\begin{eqnarray} M_{f,2}(\alpha,\beta) =I^{D_1}_{1,1}(\alpha,\beta)+I^{D_2}_{1,1}(\alpha,\beta)+O\left( T^{\frac12+\theta+\varepsilon}\right). \label{eqn:momentintegralterND} \end{eqnarray}

\subsection{Diagonal contribution}
We begin with a useful lemma.  
\begin{lemma} \label{lem:DLlaurent} The Laurent series of the meromorphic function $s\mapsto \frac{L(f\times f,s)}{\zeta^{(N)}(2s)}$ about $s=1$ can be written as
$$ \frac{L(f\times f,1+s)}{\zeta^{(N)}(2(1+s))} =\frac{\a_f/2}{s}+\b_f/2+O(s).$$
\end{lemma}
Proof: We have\begin{small}
\begin{eqnarray*}
\frac{L(f\times f,1+s)}{\zeta^{(N)}(2(1+s))} = \frac{1}{\prod_{p|N}\left( 1+\frac{1}{p^{1+s}}\right)}\frac{\zeta(1+s)}{\zeta(2(1+s))}L(Sym^2f,1+s) \\ = \frac{N}{\nu(N)}\left(1+s\sum_{p|N}\frac{\ln p}{p+1} +O(s^2)\right)\left( L(Sym^2f,1)+sL'(Sym^2f,1)+O(s^2)\right)\frac{\left( 1-2\frac{\zeta'(2)}{\zeta(2)}+O(s^2)\right)}{\zeta(2)}\\ \times\left(\frac{1}{s}+\gamma +O(s) \right).
\end{eqnarray*}
\end{small}An easy calculation gives the result. 
\begin{flushright}
$\Box$
\end{flushright}

\begin{lemma}  \label{lem:expressionavecresidu} Let $\alpha, \beta \ll \ln T$ be complex numbers. We have
$$I^{D_1}_{1,1}(\alpha,\beta) = \int_\R w(t) \underset{s=0}{Res \,}\left[\frac{G(s)}{s}g_{\alpha,\beta}(s,t) \frac{L(f\times f, 1+\alpha+\beta+2s)}{\zeta^{(N)}(2(1+\alpha+\beta+2s))} \right]dt +O\left( T^{1/2}\right).$$
\end{lemma}
Proof: By the definition of $I^{D_1}_{1,1}(\alpha,\beta)$, we can write
\begin{small}
\begin{eqnarray*}
I^{D_1}_{1,1}(\alpha,\beta) & = & \sum_{m=n}\frac{\lambda_f(m)\lambda_f(n)}{m^{\frac12+\alpha}n^{\frac12+\beta}}\int_\R w(t)V_{\alpha,\beta}(mn,t)dt \\ 
                            & = & \int_\R w(t)\frac{1}{2i\pi} \int_{(\sigma)}\frac{G(s)}{s}g_{\alpha,\beta}(s,t) \sum_{m\geq 1} \frac{\lambda_f(m)^2}{m^{1+\alpha+\beta +2s}}ds dt.                           
\end{eqnarray*}\end{small}Since $L(f\times f,s)=\zeta^{(n)}(2s)\sum_{m\geq 1}\frac{\lambda_f(m)^2}{m^{s}}$, then for all positive real number  $\sigma$ we get
\begin{small}$$I^{D_1}_{1,1}(\alpha,\beta) = \int_\R w(t)\frac{1}{2i\pi} \int_{(\sigma)}\frac{G(s)}{s}g_{\alpha,\beta}(s,t) \frac{L(f\times f, 1+\alpha+\beta+2s)}{\zeta^{(N)}(2(1+\alpha+\beta+2s))}ds dt.$$\end{small}If $\alpha+\beta \neq 0$, we specialise  \begin{small}$$G(s)= e^{s^2}\frac{(\alpha+\beta)^2-(2s)^2}{(\alpha+\beta)^2}.$$\end{small}in order to insure $G\left(- \frac{\alpha+\beta}{2}\right)=0$. We move the integration line from $\Re s=\sigma$ to $\Re s = -A$, with $A= \frac{1}{4}+\frac{\alpha+\beta}{2}$, crossing a pole at $s=0$. Thus,
\begin{small}$$I^{D_1}_{1,1}(\alpha,\beta) = \int_\R w(t) \underset{s=0}{Res \,}\left[\frac{G(s)}{s}g_{\alpha,\beta}(s,t) \frac{L(f\times f, 1+\alpha+\beta+2s)}{\zeta^{(N)}(2(1+\alpha+\beta+2s))} \right]dt +O\left( T^{1-2A}\right).$$\end{small}
\begin{flushright}
$\Box$
\end{flushright}
In order to calculate the residue at $s=0$, which appears in the previous lemma, we split our proof according to the multiplicity of this pole.
\paragraph{Double pole case}
In this paragraph, we assume $\alpha+\beta=0$. Then the pole at $s=0$ in the previous lemma has multiplicity $2$.
\begin{lemma} \label{lem:expresiducassimple} We have
\begin{small}
\begin{eqnarray*}
\underset{s=0}{Res \,}\left[\frac{G(s)}{s}g_{\alpha,\beta}(s,t) \frac{L(f\times f, 1+\alpha+\beta+2s)}{\zeta^{(N)}(2(1+\alpha+\beta+2s))} \right]  = \frac{\a_f}{2}\ln\left( \frac{t\sqrt{N}}{2\pi}\right)+\frac{\b_f}{2}
\end{eqnarray*}
\end{small}
\end{lemma}
Proof: We compute the following asymptotic behaviour at $s=0$.
\begin{small}
\begin{eqnarray*}G(s)\left( \frac{t\sqrt{N}}{2\pi}\right)^{2s} \frac{L(f\times f,1+2s)}{\zeta^{(N)}(2(1+2s))}  & = &  \left[1+O(s^2)\right]\left[1+2s\ln \left( \frac{t\sqrt{N}}{2\pi}\right)+O(s^2) \right]\left[\frac{\a_f/2}{2s}+\b_f/2+O(s)\right] \\
& = & \frac{\a_f/2}{s} +  \frac{\a_f}{2}\ln\left( \frac{t\sqrt{N}}{2\pi}\right)+\frac{\b_f}{2} +O(s).
\end{eqnarray*}
\end{small}
\begin{flushright}
$\Box$
\end{flushright} These results prove theorem \ref{th:moment2nonramolli} when $\alpha+\beta=0$. Precisely, by relation (\ref{eqn:momentintegralterND}) and since  $I^{D_1}_{1,1}(\alpha,\beta) = I^{D_2}_{1,1}(\alpha,\beta)$ in the considered case,   the following corollary comes from lemmas \ref{lem:expressionavecresidu} and \ref{lem:expresiducassimple}. 
\begin{corollary} Let $\alpha \ll L^{-1}$ be a complex number. For all $\varepsilon >0$, we have
$$ M_{f,2}(\alpha,-\alpha) = a_f \int_{\R}w(t)\ln tdt +b_f \int_\R w(t) dt +O\left( T^{\frac12+\theta+\varepsilon}\right).$$
\end{corollary}

\paragraph{Simple pole case}

In this paragraph, we assume $\alpha+\beta\neq 0$ and $G(s)= e^{s^2}\frac{(\alpha+\beta)^2-(2s)^2}{(\alpha+\beta)^2}$. We also set $w_2(t)=w(t)\left( \frac{t\sqrt{N}}{2\pi}\right)^{-2(\alpha+\beta)}$.
\begin{lemma} Let $\alpha, \beta \ll L^{-1}$ be complex numbers. For all $\varepsilon >0$, we have
\begin{small}
$$M_{f,2}(\alpha,\beta) = \frac{L(f\times f, 1+\alpha+\beta)}{\zeta^{(N)}(2(1+\alpha+\beta))} \widehat{w}(0)+ \frac{L(f\times f, 1-\alpha-\beta)}{\zeta^{(N)}(2(1-\alpha-\beta))} \widehat{w}_2(0)+O\left( T^{\frac12+\theta+\varepsilon}\right).$$
\end{small}
\end{lemma}
Proof: Since the pole at $s=0$, which appears in  lemma \ref{lem:expressionavecresidu}, is simple, we have
 \begin{small}
\begin{eqnarray*}
\underset{s=0}{Res \,}\left[\frac{G(s)}{s}g_{\alpha,\beta}(s,t) \frac{L(f\times f, 1+\alpha+\beta+2s)}{\zeta^{(N)}(2(1+\alpha+\beta+2s))} \right]  = \frac{L(f\times f, 1+\alpha+\beta)}{\zeta^{(N)}(2(1+\alpha+\beta))}.
\end{eqnarray*}
\end{small}Then, by lemma \ref{lem:expressionavecresidu}
\begin{small}
$$I^{D_1}_{1,1}(\alpha,\beta) = \frac{L(f\times f, 1+\alpha+\beta)}{\zeta^{(N)}(2(1+\alpha+\beta))} \widehat{w}(0)+O\left( T^{1/2}\right) .$$
\end{small}Up to changing $w$ by $w_2$, we have $I^{D_2}_{1,1}(\alpha,\beta) = I^{D_1}_{1,1}(-\beta,-\alpha)$, thus
\begin{small}
$$I^{D_2}_{1,1}(\alpha,\beta) = \frac{L(f\times f, 1-\alpha-\beta)}{\zeta^{(N)}(2(1-\alpha-\beta))} \widehat{w}_2(0)+O\left( T^{1/2}\right).$$
\end{small} \begin{flushright}
$\Box$
\end{flushright}
As a result, we obtain theorem \ref{th:moment2nonramolli} when $\alpha+\beta\neq 0$. 
\begin{corollary}  Let $\alpha, \beta \ll L^{-1}$ be  complex numbers. For all $\varepsilon >0$, we have
$$M_{f,2}(\alpha,\beta) = \a_f \int_{-\infty}^{+\infty} w(t)\ln t dt + \left[\b_f+\a_f \ln\left( \frac{\sqrt{N}}{2\pi}\right)\right] \widehat{w}(0) \\ +O\left(|\alpha+\beta|T(\ln T)^2+T^{\frac12+\theta+\varepsilon} \right). $$
\end{corollary}
Proof:   Thanks to lemma \ref{lem:DLlaurent}, the previous lemma gives
\begin{small}
\begin{eqnarray*}
M_{f,2}(\alpha,\beta)= \frac{\a_f/2}{\alpha +\beta} \left(\widehat{w}(0)-\widehat{w}_2(0) \right)+\frac{\b_f}{2}\left(\widehat{w}(0)+\widehat{w}_2(0) \right)+O\left(T|\alpha+\beta|+T^{\frac12+\theta+\varepsilon} \right).
\end{eqnarray*} \end{small} Furthermore,  
\begin{small}
\begin{eqnarray*}
\widehat{w}_2(0) &  = & \int_{\R}w(t)\left(\frac{t\sqrt{N}}{2\pi} \right)^{-2(\alpha+\beta)}dt = \int_{\R}w(t)\left(1-2(\alpha+\beta)\ln \left(\frac{t\sqrt{N}}{2\pi} \right)+O\left(|\alpha+\beta|^2\ln^2t \right)\right) \\
    & =& \widehat{w}(0)-2(\alpha+\beta)\int_\R w(t)\ln \left(\frac{t\sqrt{N}}{2\pi} \right)dt +O\left(|\alpha+\beta|^2T\ln^2T \right).
\end{eqnarray*}
\end{small}
An easy calculation gives the result.
\begin{flushright}
$\Box$
\end{flushright}
We remark that, up to changing $\Delta=T/\ln T$ by $T/(\ln T)^2$, we obtain corollary \ref{cor:corZhang}.

\subsection{Conjecture of Conrey, Farmer, Keating, Rubinstein and Snaith} \label{ss:conjCFKRS}
We can find in \cite{CFKRS} numerous conjectures related to integral moments of $L$-functions. In particular, conjecture 2.5.4 predicts the asymptotic behaviour of any even integral moment of a primitive $L$-function on the critical line. 
\begin{conjecture} \label{conj:conj1mi} Let $\mathcal{L}(s)$ be a primitive $L$-function. Let $k$ be a positive integer. Then for any ``suitable'' weight function $g$, we have
$$\int_{-\infty}^{+\infty} \left|\mathcal{L}\left(\frac{1}{2}+it \right)\right|^{2k}g(t) dt = \int_{-\infty}^{+\infty} P_k \left(w \ln \left(\frac{Q^{2/w}t}{2} \right) \right)(1+O(t^{-1/2+\varepsilon}))g(t)dt. $$
where $w$  and $Q$ are respectively  the degree and the conductor of $\mathcal{L}$ and $P_k$ is an explicit polynomial of degree $k^2$.
\end{conjecture}
In this paper, we consider $L$-functions of  holomorphic primitive cusp forms  of even weight, square-free level $N$  and trivial character. The degree of such an $L$-function is $w=2$ and the conductor is $Q=\sqrt{N}/\pi$ (confer equation (\ref{eqn:facteurlocalinfiniFM})). The following conjecture is a simple rewriting  of  conjecture \ref{conj:conj1mi} in this case for the second integral moment and when $g(t)=r(t/T)$ with $r$   a smooth function compactly supported in $[1,2]$.     
\begin{conjecture} \label{conj:conj2mi} Let $f$ be a holomorphic primitive cusp forms  of even weight, square-free level $N$  and trivial character. Then, for any $\varepsilon >0$, we have
 $$\int_{-\infty}^{+\infty} \left|L\left(f,\frac{1}{2}+it \right)\right|^{2}g(t) dt = \int_{-\infty}^{+\infty} P_1 \left(2 \ln \left(\frac{t\sqrt{N}}{2\pi} \right) \right) g(t)dt+O\left(T^{\frac12+\varepsilon}\right) $$
 with
 $$P_1(x)=\frac{-1}{(2i\pi)^2}\oint_{|z_1|=r_1} \oint_{|z_2|=r_2} \frac{L(f\times f, 1+z_1-z_2)}{\zeta^{(N)}(2(1+z_1-z_2))}\frac{(z_2-z_1)^2}{z_1^2z_2^2}e^{\frac{x}{2}(z_1-z_2)}dz_1dz_2,  $$
 for any small positive real numbers $r_1$ and $r_2$ (ie $r_1+r_2 < 1$).
\end{conjecture}
In order to compare our corollary \ref{cor:moment2nonramollidc} with this conjecture, we have to compute $P_1$. Choosing $r_1 \neq r_2$, and since $\frac{(z_2-z_1)^2}{z_1^2z_2^2} = \frac{1}{z_1^2}-\frac{2}{z_1z_2}+\frac{1}{z_2^2}$, we have
\begin{eqnarray*} 
P_1(x) & = & \frac{2}{(2i\pi)^2}\oint_{|z_1|=r_1} \oint_{|z_2|=r_2} \frac{L(f\times f, 1+z_1-z_2)}{\zeta^{(N)}(2(1+z_1-z_2))}\frac{1}{z_1z_2}e^{\frac{x}{2}(z_1-z_2)}dz_1dz_2 \\ 
 & = & \frac{2}{2i\pi}\oint_{|z_2|=r_2}\frac{L(f\times f, 1-z_2)}{\zeta^{(N)}(2(1-z_2))}e^{-\frac{x}{2}z_2}\frac{dz_2}{z_2}
\end{eqnarray*}
Moreover, since 
\begin{small}$$\frac{L(f\times f, 1-z_2)}{\zeta^{(N)}(2(1-z_2))}e^{-\frac{x}{2}z_2} = \left( \frac{-\a_f/2}{z_2}+\b_f +O(z_2^2) \right)\left( 1-\frac{x}{2}z_2+O(z_2^2)\right) =\frac{-\a_f/2}{z_2} +\left( \frac{\b_f}{2}+\frac{x\a_f}{4}\right) +O(z_2^2),$$\end{small}we obtain
$$P_1(x)=\frac{\a_f}{2}x+\b_f.$$
To conclude, we can see that the main terms are similar in corollary \ref{cor:moment2nonramollidc} and conjecture \ref{conj:conj2mi} and, assuming the Ramanujuan-Petersson conjecture, the error terms are also equal. 
\nocite{HY}

\newcommand{\etalchar}[1]{$^{#1}$}

\noindent Damien BERNARD \\
Université Blaise Pascal \\
Laboratoire de Mathématiques \\
Campus des Cézeaux \\
BP 80026\\
63171 Aubière cedex, France\\
E-mail address: damien.bernard@math.univ-bpclermont.fr

\end{document}